\documentclass[12pt]{article}
\usepackage{amsmath,amssymb,fullpage}
\usepackage[applemac]{inputenc}

\newenvironment{myenumerate}{%

\begin{enumerate}}{\end{enumerate}}

\newcommand{\dproof}{\noindent {Proof.} \quad}
\newcommand{\fproof}{\hfill $\square$ \bigskip}

\newtheorem{example}{Example}[section]
\newtheorem{theorem}[example]{Theorem}

\newtheorem{problem}[example]{Problem}
\newtheorem{remark}[example]{Remark}

\numberwithin{equation}{section}

\def\RB{\mathbb{R}}

\def\FB{\mathbb{F}}

\def\BB{\mathbb{B}}
\def\GB{\mathbb {G}}

\def\FC{\mathcal{F}}

\def\EC{\mathcal{E}}
\def\AC{\mathcal{A}}

\def\VC{\mathcal V}

\def\XC{\mathcal X}

\def\GC{\mathcal G}

\def\ess{\mathop{ess \;sup}}
\def\esi{\mathop{ess \;inf}}

\def\1B{\text{1\!\!I}}

\def\tX{\tilde{X}}
\def\tsi{\tilde{\sigma}}

\def\hX{\hat{X}}

\def\hY{\hat{Y}}
\def\hu{\hat{u}}

\def\hb{\hat{b}}
\def\hh{\hat{\cal H}}

\def\hH{\hat{H}}
\def\hf{\hat{f}}

\def\hp{\hat{p}}
\def\hth{\hat{\theta}}
\def\hq{\hat{q}}
\def\hr{\hat{r}}
\def\hw{\hat{w}}

\def\hxi{\hat{\xi}}
\def\hla{\hat{\lambda}}
\def\hsi{\hat{\sigma}}

\def\EE{E}

\def\cA{{\mathcal A}}

\def\be{{\beta}}

\def\la{{\lambda}}
\def\si{{\sigma}}

\def\al{{\alpha}}

\def\be{{\beta}}

\def\si{{\sigma}}

\def\la{{\lambda}}

\def \eref#1{\hbox{(\ref{#1})}}

\def\si{{\sigma}}
\def\al{{\alpha}}

\def\hX{\hat{X}}

\def\hY{\hat{Y}}

\def\hu{\hat{u}}
\def\hw{\hat{w}}
\def\hb{\hat{b}}
\def\hf{\hat{f}}
\def\hh{\hat{h}}
\def\hH{\hat{\cal H}}

\def\hH{\hat{H}}
\def\hp{\hat{p}}
\def\hth{\hat{\theta}}
\def\hq{\hat{q}}
\def\hr{\hat{r}}
\def\hxi{\hat{\xi}}
\def\hla{\hat{\lambda}}
\def\hsi{\hat{\sigma}}

\def\hG{\hat{G}}

\def\ess{\mathop{ess \;sup}}
\def\esi{\mathop{ess \;inf}}

\def\tX{\tilde{X}}
\def\tsi{\tilde{\sigma}}

\def\RB{\mathbb{R}}

\def\FB{\mathbb{F}}

\def\BB{\mathbb{B}}
\def\GB{\mathbb{G}}

\def\GC{\mathcal{G}}
\def\FC{\mathcal{F}}

\def\EC{\mathcal{E}}
\def\AC{\mathcal{A}}

\def\VC{\mathcal V}
\def\XC{\mathcal X}

\let\Section=\section
\def\section{\setcounter{equation}{0}\Section}

\begin{document}

\title{Singular mean-field control games with applications to optimal harvesting and investment problems}

\date{ 6 June 2014}

\author{
Yaozhong Hu$^{1,2}$, Bernt \O ksendal$^{2,3}$, Agn\`es Sulem$^{4,5,2}$}

\footnotetext[1]{Department of Mathematics, University of Kansas, Lawrence, Kansas 66045, USA.  
Email: {\tt yhu@ku.edu}. \\
Y. Hu is partially supported by a grant from the Simons Foundation \#209206
and by a General Research Fund of University of Kansas.}
\footnotetext[2]{Department of Mathematics, University of Oslo, P.O. Box 1053 Blindern, N--0316 Oslo, Norway.\\
  Email: {\tt oksendal@math.uio.no}. \\
The research leading to these results has received funding from the European Research Council under the European Community's Seventh Framework Programme (FP7/2007-2013) / ERC grant agreement no [228087].}
\footnotetext[3]{ 
   Norwegian School of Economics (NHH), Helleveien 30, N--5045 Bergen, Norway.}  
\footnotetext[4]{ INRIA Paris-Rocquencourt, Domaine de Voluceau, BP 105, Le Chesnay Cedex, 78153, France;   \\
 Email: {\tt agnes.sulem@inria.fr}}

\footnotetext[5]{
Universit\'e Paris-Est, F-77455 Marne-la-Vallée, France.}

\maketitle



\begin{abstract}This paper studies singular mean field control problems
and singular mean field stochastic differential games.
Both sufficient  and necessary conditions for the optimal controls and
for the Nash equilibrium are obtained. Under some  assumptions the optimality conditions
for singular mean-field control are
reduced to a reflected Skorohod problem, whose solution is proved to exist
uniquely. Some examples are given. In particular, a simple singular mean-field
 game is studied where the Nash equilibrium exists but is not unique.
\end{abstract}

\maketitle

 \section{Introduction}\label{sec1}

The irreversible investment problem is a classical problem in economics,
with a long history. In short, the problem is the following.
A factory is facing an increased demand for its product. Should it invest in more production
capacity to meet the demand? The problem is that buying additional production capacity is an
 expensive, irreversible investment (usually production equipment cannot easily be sold
 after use) and the future demand for the product is uncertain. So the risk is that the
 factory ends up having paid for an additional capacity it does not need. On the other hand,
 if the factory does not increase the capacity, it might  miss an opportunity for an increased sale.

This is a classical problem that has been studied by many authors in different contexts.
See e.g. Pindyck (1988, 1991, 1991),\cite{P1,P2, P3} , Kobila (1993) \cite{K} and the references therein.\\
Mathematically the problem can be formulated as a singular control problem. In this paper we
study such singular control problems in the context of \emph{mean-field It\^o processes}
and under \emph{model uncertainty}. We interpret model uncertainty in the sense of Knight uncertainty,
i.e. uncertainty about the underlying probability measure. Using the Girsanov theorem for It\^o
processes we can parametrize the family of densities of possible underlying probability measures by a stochastic
process $\theta(t)$. This leads to a stochastic differential game in which one of the players is the
investor controlling the investment strategy and the other player controls the model by choosing the
model parameter $\theta$. Since the investment is irreversible, the control of the investor is a
\emph{singular} control, i.e. a non-decreasing non-negative stochastic process $ \xi(t)$. Thus we
are dealing with a \emph{singular mean-field stochastic differential game}. 

Recently,  there have been several papers dealing with mean-field control problems. See e.g.  Meyer-Brandis et al (2012) \cite{MOZ},  Anderson \& Djehiche (2011) \cite{AD} and Hamad\`ene (1998) \cite{H}.
A recent paper dealing with mean-field \emph{singular} control is L. Zhang (2012) \cite{Z}.
Our paper extends this paper in several directions:
For example, we consider more general mean-field operators. And we allow the profit rate
$f $ (see below) to depend on the state $X$, on the mean-field term $Y$,   and on the singular
control $\xi$. We also allow both the coefficient $\lambda$ in the singular part of the state
equation and the singular cost coefficient $h$ in the performance functional to depend on
the state  $X$.
Moreover,  we consider general games between two players of such singular control problems with
asymmetric information.

Our paper is organized as follows:
 In Section 2 we present three motivating examples, 
 In Section 3 we formulate the general mean-field singular stochastic control problem and prove a
sufficient maximum principle and a necessary maximum principle.
 In Section 4 we reduce the maximum principle to a Skorohod problem and prove the existence and
uniqueness of the solution.
 In Section 5  we prove a \emph{sufficient} maximum principle for general singular mean-field
  stochastic games, and we obtain as a  corollary a
 corresponding maximum principle for zero-sum games.
 In Section 6 we apply the results above to an optimal harvesting problem of a mean-field system and
 to model uncertainty singular control, in particular model uncertainty irreversible investment type problems.

\section{Three motivating examples}
\subsection{Optimal harvesting from a mean-field system}

Suppose we model the size $X^0(t)  $ of an unharvested population at time t by an equation of the form
\begin{align}\label{eq2.1a}
dX^{0}(t) &= E[X^{0}(t)] b(t)dt + X^{0}(t)\sigma(t) dB(t)
; \; \; t \in [0,T] \nonumber\\
X^0(0) &= x > 0.
\end{align}
Here, and in the following,  $B(t) = B(t,\omega)$
 is a Brownian motion
on a filtered probability space $(\Omega,   \FC, \mathbb{F}:=\{\FC_t\}_{t \geq 0}, P)$ satisfying the usual conditions.  $P$ is a reference probability measure.
We assume that $b(t)$ and $\sigma(t)$ are given predictable processes.
We may regard \eqref{eq2.1a} as a limit as $n \rightarrow \infty$ of a large population interacting system of the form
\begin{equation}
dx^{i,n}(t)= [\frac{1}{n} \sum_{j=1}^n x^{j,n}(t)] b(t) dt
+ x^{i,n}(t)\sigma(t) dB^i(t) , \quad i = 1,2, ..., n.
\end{equation}
Thus the mean-field term $E[X(t)]$ represents an approximation to the weighted average $ \frac{1}{n}\sum_{j=1}^n  x^{j,n}(t)$ for large $n$.
Now suppose we introduce harvesting of the population. The size of the harvested population $X(t) = X^\xi(t)$ at time $t$ can then be modeled by a mean-field singular control stochastic differential equation of the form
\begin{align}\label{eq2.1}
dX(t) &= E[X(t)] b(t)dt + X(t )\sigma(t) dB(t)
- \lambda_0(t) d\xi(t); \; \; t \in [0,T] \nonumber\\
X(0 ) &= x > 0
\end{align}
where $\xi(t)$ is a non-decreasing predictable process with $\xi(0^-)=0$, representing the harvesting effort, while $\lambda_0(t) > 0$ is a given harvesting efficiency coefficient.

The \emph{performance functional} is assumed to be of the form
\begin{equation}\label{eq2.2}
J(\xi) = E\left[\int_o^T h_0(t)X(t) d\xi(t) + KX(T)\right],
\end{equation}
where $h_0(t)$ is a given adapted price process and $K = K(\omega)$ is
a given salvage price, assumed to be $\mathcal{F}_T$-measurable.
The problem is to find $\xi^*$ such that
\begin{equation} \label{eq2.3}
J(\xi^*) = \sup_{\xi} J(\xi)\,.
\end{equation}
Such a  process $\xi^*$ is called an optimal singular control.
This is an example of a \emph{mean-field singular control problem}.
We will return to this problem in Section 6.1.

\subsection{Optimal irreversible investments under model uncertainty}

Let $\xi(t)$ denote the production rate capacity of a production plant and let
$D(t)$ denote the demand rate at time $t$. At any time $t$ the production
capacity can be increased by $d\xi(t)$ at the price $\lambda_0(t,D(t))$ per
capacity unit. The number of units sold per time unit is the minimum of the demand D(t) and the capacity $\xi(t)$.The total expected net profit of the production is assumed to be
\begin{align}
J(\xi,\theta)
=& E^{Q_\theta}\Bigg[ \int_0^T a(t,E[\varphi(D(t))]) \min [D(t), \xi(t)] dt \nonumber\\
&\qquad + g(D(T)) - \int_0^T \lambda_0 (t, D(t)) d\xi(t) \Bigg]\,, \label{eq2.4}
\end{align}
where $g(D(T))$ is some salvage value of the closed-down production plant, $\varphi$ is a given real
function and $a(t,E[\varphi(D(t))])$ is the unit sales price of the production. Here $\{Q^\theta\}_{\theta \in \Theta}$
is a family of probability measures representing the model uncertainty. We let $\AC_\mathbb{G}$
denote the set of right-continuous, non-decreasing $\mathbb{G}$-adapted processes $\xi(\cdot)$ with
$\xi(0^-) = 0$, where $\mathbb{G}:= \{ \mathcal{G}_t \}_{t\geq0}$ is a given subfiltration of $\mathbb{F}$, in the sense that $\GC_t \subseteq \FC_t$ for all $t$. Heuristically, $\mathcal{G}_t$ represents the
information available to the investor at time $t$.
We assume that the demand $D(t)$ is given by a jump diffusion of the form
\begin{equation}\label{eq2.5}
\begin{cases}
dD(t) & = \displaystyle D(t^-) \left[ \alpha (t,\omega) dt + \beta(t,\omega)dB(t)
\right], 0 \leq t \leq T\\
D(0) & > 0
\end{cases}
\end{equation}
where $\alpha(t,\omega), \beta(t,\omega)$ are given $\mathbb{F}$-adapted processes.
We want to maximize the expected total net profit under the worst possible scenario, i.e. find $(\xi^*, \theta^*) \in \AC_\GB \times \Theta$ such that
\begin{equation}\label{eq2.6}
\sup_{\xi \in \AC_\GB} \left\{ \inf_{\theta \in \Theta} J(\xi, \theta)\right\} = \inf_{\theta \in \Theta} \left\{ \sup_{\xi \in \AC_\GB} J(\xi, \theta)\right\} = J (\xi^*, \theta^*).
\end{equation}

This is an example of a (partial information) {\it singular control game} of a jump diffusion. Note that the system is non-Markovian, both because of the mean-field term
and the partial information constraint. See Sections 6.2-6.4.

\subsection{A mean field singular game}\label{ex2.3}
Suppose the demand $X(t)$ for a certain product at time $t$ is given by a mean field
SDE of the form
\[
\left\{
\begin{array}{lll}
dX(t)&=&\EE [X(t)] b(t) dt+X(t)\si(t) dB(t)\\
X(0)&=&x>0\,.
\end{array}\right.
\]
There are two competing companies producing this product, with production rate
capacities represented by nondecreasing adapted processes $\xi_1$,  $\xi_2$,
respectively. The expected profit of the company $i$ is assumed to have to the form
\[
J_i(\xi_1, \xi_2)=\EE \left[ \int_0^T \pi(t) \min(X(t), \xi_1(t)+ \xi_2(t))dt+\int_0^T
h_i(t)d\xi_i(t)\right]\,,
\]
where $\pi(t)>0$ is the price per unit sold and $h_i(t)<0$ the production cost per unit
for the factory $i$, $i=1,2$.
We want to find a \emph{Nash equilibrium}, i.e. a pair $(\xi^*_1, \xi^*_2)\in \cA_1\times \cA_2$
such that
\[
\sup_{\xi_1\in \cA_1} J_1(\xi_1, \xi_2^*)=J_1(\xi_1^*, \xi_2^*)
\]
and
\[
\sup_{\xi_2\in \cA_2} J_2(\xi_1^*, \xi_2 )=J_2(\xi_1^*, \xi_2^*)\,,
\]
where $\mathcal{A}_i$ is the family of admissible controls $\xi_i$ for company number $i; i=1,2$.
We will return to this problem in Section 6.5.

\section{Maximum principle for singular mean field control problems}\label{sec3}

\subsection{Problem statement}
We first recall some basic concepts and results from Banach space theory.
 Let $V$ be an open subset of a Banach space $\XC$ with norm $\| \cdot \|$ and let $F : V \rightarrow \RB$.
\begin{myenumerate}
\item We say that $F$ has a directional derivative (or G\^ateaux derivative) at $x \in \XC$ in the direction $y \in \XC$ if
$$D_y F(x) := \lim_{\varepsilon \rightarrow 0} \frac{1}{\varepsilon} (F(x + \varepsilon y) - F(x))$$ exists.
\item We say that $F$ is  Fr\'echet differentiable at $x \in V$ if there exists a linear map
$L: \XC \rightarrow \RB$
such that
$$\lim_{\substack{h \rightarrow 0 \\ h \in \XC}} \frac{1}{ \|h\|} | F(x+h) - F(x) - L(h)| = 0.$$
In this case we call $L$ the {\it gradient} (or Fr\'echet derivative) of $F$ at $x$ and we write
$$L = \nabla_x F.$$
\item If $F$ is Fr\'echet differentiable, then $F$ has a directional derivative in all directions $y \in \XC$ and
$$D_y F(x) = \nabla_x F(y) =: \langle \nabla_xF,y\rangle .$$
\end{myenumerate}
In particular, if $\XC=L^2(P)$ the Fr\'echet derivative of $F$ at $X \in L^2(P)$,
denoted by $\nabla_X F$, is a bounded linear functional on $L^2(P)$, which we can identify with a random variable in $L^2(P)$.
For example, if $F(X) = E[\varphi(X)]; X \in L^2(P)$, where $\varphi$ is a real  $C^1$- function such that $\varphi(X) \in L^2(P)$ and $  \frac{\partial \varphi}{\partial x} (X) \in L^2(P)$, then
$\nabla_X F = \frac{\partial \varphi}{\partial x} (X)$ and $\nabla_XF(Y) =
\langle \frac{\partial \varphi}{\partial x}(X),Y\rangle =E[\frac{\partial \varphi}{\partial x}(X)Y]$ for $Y \in L^2(P).$
\vskip 0.3 cm

Consider a mixed regular and singular controlled  system with state process $X(t) = X^{\xi, u}(t)$ of the form
\begin{align}\label{eq3.1a}
dX(t) & = b(t,X(t),Y(t),\xi(t), u(t),\omega)dt + \sigma(t,X(t),Y(t),\xi(t), u(t),\omega)dB(t) \nonumber \\
&
 + \lambda(t,X(t),u(t),\omega)d\xi(t),
\end{align}
where
\begin{equation}\label{eq3.2}
Y(t) = F(X(t, \cdot))
\end{equation}
and $F$ is a Fr\'echet differentiable operator on $L^2(P)$.
We assume that all the coefficients $b, \sigma, \lambda, f, g$ and $h$ are Fr\'echet differentiable ($C^1$) with respect to $x, y, \xi,u$ with derivatives in $L^2(m \times P)$, where $m$ denoted Lebesgue measure on $[0,T]$.
Note that we allow the coefficients to depend on both controls $\xi$ and $u$. This might be relevant, for example, in harvesting models. See \eqref{ex4.3}.
The \emph{performance functional}   is assumed to be of the form
\begin{align}\label{eq3.3}
J (\xi,u) & = E \left[ \int_0^T f (t,X(t),Y(t),\xi(t), u(t), \omega)dt + g (X(T),Y(T),\omega)\right. \nonumber \\
 &\qquad  + \left. \int_0^T h (t,X(t),u(t),\omega) \xi(dt)\right] \,.
 \end{align}
 We may interpret the function $f $ as a profit rate, $g $ as a bequest or salvage value function and $h $
 as a cost rate for the use of the singular control $\xi$.

  We want to find  $(\xi^*, u^*) \in \AC $   such that
 \begin{equation}\label{eq3.3b}
 J (\xi^*,u^*)=\sup_{(\xi,u)  \in \AC } J (\xi,u)\,.
 \end{equation}
Here $\AC$ is a given family of $\GB $-predictable processes such that the corresponding state equation has a
unique solution $X$ such that $\omega \rightarrow X(t,\omega) \in L^2(P)$ for all t.
We let $A $ denote the set of possible values of $u(t); t \in [0,T] $ when $(\xi,u) \in \AC.$

\subsection{A sufficient maximum principle for singular mean field control}

 In this subsection we prove a sufficient maximum principle for the singular control games described above.
 To this end, define the  {\it Hamiltonians} $H$ as follows:
\begin{align}\label{eq5.2a}
 H& (t,x,y,\xi,u,p,q)(dt, \xi(dt)) \nonumber \\
 & = H_{0} (t,x,y,\xi,u,p,q) dt + \{ \lambda(t,x,u) p + h(t,x,u)\} \xi(dt) \,,
 \end{align}
 where
 \begin{align}\label{eq5.3a}
 H_0(t,x,y,w)= f(t,x,y,\xi, u) + b(t,x,y,\xi, u)p + \sigma(t,x,y,\xi, u)q\,.
 \end{align}
The associated mean-field BSDE for the adjoint processes is
\begin{equation}\label{eq3.4b}
 \begin{cases}
 dp(t) & = \displaystyle - \frac{\partial H_{0}}{\partial x} (t,X(t),Y(t),\xi(t), u(t) ,p(t),q(t))dt\\
 &\displaystyle- \frac{\partial H_{0}}{\partial y} (t,X(t), Y(t), \xi(t), u(t) ,p(t),q(t))\nabla_{X(t)} F)dt
  + q(t) dB(t) \\
  p(T) & = \displaystyle \frac{\partial g}{\partial x}(X(T),Y(T))+E[\frac{\partial g}{\partial y}(X(T),Y(T))]\nabla_{X(T)} F.
\end{cases}
\end{equation}

The  sufficient maximum principle for this singular mean field control is  stated as follows.

\begin{theorem}[Sufficient maximum principle for mean-field singular control]\label{th5.1a}
Let $\hxi, \hu \in \AC$, with corresponding solutions $\hX(t), \hY(t), \hp(t),\hq(t)$ of \eqref{eq3.1a} and \eqref{eq3.4b}.
Suppose the following conditions hold
\begin{itemize}
\item
The function
 \begin{equation}
 X,\xi, u \rightarrow H(t,X,F(X),\xi, u,\hp(t),\hq(t))\label{e.3.8a}
 \end{equation}
 is concave with respect to $(X,\xi, u) \in L^2(P) \times \mathbb{R} \times \mathbb{R}$ for all $t\in [0, T]$.\\

\item (The conditional maximum property)
\begin{align}\label{eq4.43a}
\ess_{v \in A} E [H(t, \hX(t), \hY(t),\hxi(t),v, \hp(t), \hq(t)) \mid \GC_t] \nonumber \\
= E[H(t, \hX(t), \hY(t), \hxi(t), \hu(t), \hp(t), \hq(t))\mid \GC_t]\,.
\end{align}

\item (Variational inequality)
\begin{align}\label{eq4.44a}
\ess_{\xi} E [H(t, \hX(t), \hY(t),\xi,\hu(t), \hp(t), \hq(t)) \mid \GC_t] \nonumber \\
= E[H(t, \hX(t), \hY(t), \hxi(t), \hu(t), \hp(t), \hq(t))\mid \GC_t]
\end{align}

\end{itemize}

Then $(\hxi(t),\hu(t))$ is an optimal control for $J(\xi, u)$.
\end{theorem}
\dproof This theorem is a straightforward consequence of Theorem \ref{th3.1} below. We refer to the proof
there.
\fproof 

\subsection{A necessary maximum principle for singular mean field control}\label{sec5}

In the previous section we gave a verification theorem, stating that if a given control $(\hxi, \hu)$
satisfies \eref{e.3.8a}-\eref{eq4.44a}, then it is indeed optimal for the singular mean field control
problem. We now establish a
partial converse, implying that if a control $(\hxi, \hu)$ is optimal for the singular mean field
control problem, then
it is a conditional critical  point for the Hamiltonian.

To achieve this, we start with the setup of \cite{OS3} as follows.
 For $\xi \in \AC $,  let $\VC(\xi)$ denote the set of $\GB$-adapted processes $\eta$ of finite
 variation such that there exists $\delta = \delta(\xi) > 0$ satisfying
\begin{equation}\label{eq5.1a}
\xi + a \eta \in \AC  \text{ for all } a \in [0, \delta].
\end{equation}
Then for $\xi \in \AC $ and $\eta \in \VC(\xi)$ we have, by our smoothness assumptions on the coefficients,
\begin{align}\label{eq5.2a}
\lim_{a \rightarrow 0^+ }  &\frac{1}{a} \left( J (\xi + a \eta) - J(\xi)\right)
= E \left[ \int_0^T \left\{ \frac{\partial f}{\partial x} (t, X(t), Y(t), \xi(t), u(t) ) Z(t) \right. \right. \nonumber \\
 & + \left. \frac{\partial f}{\partial y} (t, X(t) Y(t), \xi(t), u(t) ) \langle \nabla_{X(t)} F, Z(t) \rangle \right\} dt \nonumber \\
  & \left.+ \int_0^T \frac{\partial f}{\partial \xi} (t, X(t), Y(t), \xi(t), u(t) ) \eta (dt) \right] \nonumber \\
  & + E \left[ \frac{\partial g}{\partial x} (X(T), Y(T)) Z(T) + \frac{\partial g}{\partial y}
  (X(T), Y(T)) \langle \nabla_{X(T)} F, Z(T) \rangle \right] \nonumber \\
  & + E \left[ \int_0^T \left\{ \frac{\partial h}{\partial x} (t, X(t), \xi(t), u(t) ) Z(t) \xi(dt)
  + \int_0^T h(t, X(t), \xi(t), u(t) ) \eta(dt)\right] \right.
  \end{align}
where
\begin{equation} \label{eq5.3a}
Z(t) := \lim_{a \rightarrow 0^+ } \frac{1}{a} \left( X^{(\xi + a \eta)}(t) - X^{(\xi)}(t)\right).
\end{equation}
Note that by the chain rule we have
\begin{align}\label{eq5.4a}
\lim_{a \rightarrow 0^+ }& \frac{1}{a} \left( Y^{(\xi + a \eta)}(t) - Y^{(\xi)}(t)\right) 
 = 
 \lim_{a \rightarrow 0^+} \frac{1}{a} \left( F(\left( X^{(\xi + a \eta)}(t)) - F(X^{(\xi)}(t)\right) \right) 
 = \langle \nabla_{X(t)} F, Z(t) \rangle \,.
\end{align}
Moreover,
\begin{align}\label{eq5.5a}
dZ(t) & = \left( \frac{\partial b}{\partial x}(t) Z(t) + \frac{\partial b}{\partial y}(t) \langle \nabla_{X(t)}F, Z(t) \rangle \right) dt \nonumber \\
& + \left( \frac{\partial \sigma}{\partial x}(t) Z(t) + \frac{\partial \sigma}{\partial y}(t) \rangle \nabla_{X(t)}F, Z(t) \rangle \right) dB(t) \nonumber \\
 & + \frac{\partial \lambda}{\partial x} (t) Z(t) \xi(dt) + \lambda (t) \eta(dt) \; ; \; Z(0) = 0\,.
\end{align}

\begin{theorem}\label{th5.1a}{(Necessary maximum principle for singular mean-field control problem)}
Suppose  $ (\hxi , \hu) \in \AC $ is optimal, i.e. satisfies \eqref{eq3.3b}. Then
\begin{equation}\label{eq5.7a}
E \left[ \frac{\partial H_{ 0}}{\partial u } (t, \hX(t), \hY(t), \hxi (t), u ,   \hp (t), \hq(t))_{u =\hu (t)}
\mid \GC _t \right] = 0\,.
\end{equation}
Moreover, the following variational inequalities hold.
\begin{equation}\label{eq3.18a}
\begin{cases}
\displaystyle E \left[ \frac{\partial f}{\partial \xi} (t, \hX(t), \hY(t), \hxi(t), \hu(t))
+ \lambda(t, \hX(t), \hu(t)) \hp(t) + h  (t, \hX(t), \hu(t)) \mid \GC _t \right] \leq 0 \\
\text{ for all } t\in [0, T] \quad
\text{and} \\
\displaystyle E \left[\frac{\partial f}{\partial \xi} (t, \hX(t), \hY(t), \hxi(t), \hu(t))
+ \lambda(t, \hX(t), \hu(t)) \hp(t) + h  (t, \hX(t), \hu(t)) \mid \GC
_t \right] d \hxi(t) =  0 \\
\text{ for all } t\in [0, T]\,.
\end{cases}
\end{equation}
\end{theorem}
\dproof To simplify the notation denote
\begin{eqnarray*}
\hf(t) &:=& f(t,\hX(t),\hY(t),\hxi(t))\\
 \hla(t)&:=&\lambda(t,\hX(t)), \qquad  \hh(t):=h(t,\hX(t)).
\end{eqnarray*}
We need to prove that if $\hxi \in \AC_{\GB}$ is optimal, i.e. if
\begin{equation}\label{eq5.14a}
\sup_{\xi \in \AC_{\GB}} J(\xi) = J(\hxi)
\end{equation}
then $\hxi$ satisfies the following variational inequalities:
\begin{equation}\label{eq5.15a}
E \left[ \frac{\partial \hf}{\partial \xi}(t) + \hla(t) \hp(t) + \hh(t) \mid \GC_t \right] \leq 0
\quad  \text{ for all } t \in [0,T]
\end{equation}
and
\begin{equation}\label{eq5.16a}
E \left[ \frac{\partial \hf}{\partial \xi}(t) + \hla(t) \hp(t) + \hh(t) \mid \GC_t \right]\hxi(dt) = 0
\quad \text{ for all } t \in [0,T].
\end{equation}

To this end, choose $\xi \in \AC_{\GB}$ and $\eta \in \VC(\xi)$ and compute
\begin{equation}\label{eq5.17a}
\frac{d}{da} J(\xi + a \eta) \mid_{a=0} = A_1 + A_2 + A_3 + A_4,
\end{equation}
where
\begin{eqnarray}
A_1 & =& E \left[ \int_0^T \left\{ \frac{\partial f}{\partial x} (t) Z(t)
+ \frac{\partial f}{\partial y}(t) \langle \nabla_{X(t)} F, Z(t) \rangle \right\} dt \right] \nonumber \\
A_2 & =& E \left[ \int_0^T \frac{\partial f}{\partial \xi} (t) \eta (dt) \right]\nonumber\\
A_3 & =& E \left[ \frac{\partial g}{\partial x} (X(T), Y(T)) Z(T)
+ \frac{\partial g}{\partial y} (X(T), Y(T)) \langle \nabla_{X(T)}F, Z(T) \rangle \right] \nonumber \\
A_4 & =& E \left[ \int_0^T \frac{\partial h}{\partial x} (t) Z(t) \xi(dt) + h (t) \eta (dt)\right].
\end{eqnarray}
By the definition of $H_0$ we have
\begin{align}\label{eq5.20}
A_1& = E \left[ \int_0^t Z(t) \left\{ \frac{\partial H_0}{\partial x}(t)
- \frac{\partial b}{\partial x}(t) p(t) - \frac{\partial \sigma}{\partial x}(t) q(t) \right\} dt \right. \nonumber \\
& \left. + \int_0^T \langle \nabla_{X(t)} F, Z(t) \rangle \left( \frac{\partial H_0}{\partial y}(t)
- \frac{\partial b}{\partial y}(t) p(t) - \frac{\partial \sigma}{\partial y}(t) q(t)\right)dt \right].
\end{align}
By the terminal condition of $p(T)$ (see \eqref{eq3.4b}) we have
\begin{align}\label{eq5.21a}
A_3 & = E [ \langle p(T), Z(T) \rangle ] \nonumber\\
 & = E \left[ \int_0^T p(t) dZ(t) + \int_0^T Z(t)dp(t) + \int_0^T q(t)
  \left( \frac{\partial \sigma}{\partial x}(t) Z(t) + \frac{\partial \sigma}
  {\partial y}(t) \langle \nabla_{X(t)}F, Z(t) \rangle \right) dt\right. \nonumber \\
 & = E \left[ \int_0^T p(t) \left\{ \frac{\partial b}{\partial x}(t) Z(t)
 + \frac{\partial b}{\partial y}(t) \langle \nabla_{X(t)}F, Z(t) \rangle  \right\} dt \right. \nonumber \\
 & \quad + \int_0^T p(t) \frac{\partial \lambda}{\partial x}(t) Z(t) \xi(dt)
 + \int_0^T p(t) \lambda(t) \eta(dt) \nonumber \\
  & \quad - \int_0^T \left( \frac{\partial H_0}{\partial x}(t) Z(t)
   + \frac{\partial H_0}{\partial y}(t) \langle \nabla_{X(t)} F, Z(t) \rangle\right) dt \nonumber \\
  & \quad + \int_0^T q(t) \left( \frac{\partial \sigma}{\partial x} (t) Z(t)
   + \frac{\partial \sigma}{\partial y}(t) \langle \nabla_{X(t)}F, Z(t) \rangle \right)dt.
\end{align}

Combining \eqref{eq5.17a}-\eqref{eq5.21a} we get
$$\frac{d}{da}J(\xi + a \eta)\mid_{a=0} = E \left[ \int_0^T
 \left\{ \frac{\partial f}{\partial \xi}(s) + \lambda(s) p(s) + h(s) \right\} \eta(ds)\right].$$

In particular, if we apply this to an optimal $\xi = \hxi$ for $J$ we get, for all $\eta \in \VC(\hxi)$,
\begin{align}\label{eq5.22a}
E \left[ \int_0^T \left\{ \frac{\partial \hf}{\partial \xi}(s) + \hla(t), \hp(s) + \hh(s) \right\}
\eta (ds) \right] = \frac{d}{da} J( \hxi + a \eta)_{a=0} \leq 0.
\end{align}
If we choose $\eta$ to be a pure jump process of the form
$$\eta(s) = \sum_{0 <t_ i \leq s} \alpha(t)$$
where $\alpha(t) > 0$ is $\GC_{t}$-measurable, \eqref{eq5.22a} gives
$$E \left[ \left( \frac{\partial \hf}{\partial \xi} (t) + \hla(t) \hp(t)
+ \hh(t)\right) \alpha \right] \leq 0
\text{ for all } t.$$

Since this holds for all such $\eta$, we conclude that
\begin{equation}\label{eq5.23a}
E \left[ \frac{\partial \hf}{\partial \xi}(t) + \hla(t) \hp(t) + \hh(t) \mid \GC_t \right] \leq 0
 \text{ for all }  t \in [0,T] .
\end{equation}
Finally, applying \eqref{eq5.22a} to
$$ \eta(dt) = \hxi(dt) \in \VC(\hxi)$$
and then to
$$ \eta(dt) = - \hxi(dt) \in \VC(\hxi)$$
we get, for all $t \in [0,T]$,
\begin{equation}\label{eq5.24a}
E \left[ \frac{\partial \hf}{\partial \xi}(t) + \hla(t) \hp(t) + \hh(t) \mid \GC_t\right]  \hxi(dt) = 0 \text{ for all } t \in [0,T].
\end{equation}
With \eqref{eq5.23a} and \eqref{eq5.24a} the proof is complete.
\fproof 

\section{The optimality conditions}
Since there have already been studies (see e.g. \cite{MOZ} and
references therein)  on the usual (nonsingular) mean field control problems,
  let us consider only the singular control $\xi$, i.e. with no regular control $u$.  The system that
we shall deal with, is described by
\begin{eqnarray}
dX(t)
&=& b(t, X(t), Y(t), \xi(t) ) dt+ \si (t, X(t), Y(t), \xi(t)  ) dB(t)\nonumber\\
&&\qquad +\la  (t, X(t), Y(t)) d\xi(t)\,.\label{e.system}
\end{eqnarray}
The performance functional is
\begin{equation}
J(\xi)=
 \EE \left[ \int_0^T  f(t, X(t), Y(t), \xi(t) ) dt +g(X(T), Y(T))+\int_0^T h(t, X(t), Y(t)) d\xi(t)\right]\,.
\end{equation}
The auxiliary backward stochastic differential equation with mean
field is
\begin{eqnarray}
dp(t)
&=& -\bigg[ f_x(t, X(t), Y(t),   \xi(t)) +b_x(t, X(t), Y(t), \xi(t)   ) p(t) \nonumber\\
&&\qquad +\si_x
 (t, X(t), Y(t),   \xi(t)) q(t)\bigg] dt
 -\bigg[ f_y(t, X(t), Y(t),   \xi(t)) \nonumber\\
&&\qquad +b_y(t, X(t), Y(t) , \xi(t) ) p(t)   +\si_y
 (t, X(t), Y(t),   \xi(t)) q(t)\bigg] \nabla_{X(t)} Fdt \nonumber\\
 &&\qquad +q(t) dB(t)\,,\label{e.bsde}\\
 p(T)&=& g_x(X(T), Y(T))+\EE \left[ g_y(X(T), Y(T)) \right] \nabla_{X(T)} F\,.
 \label{e.bsde-t}
\end{eqnarray}

To solve the above BSDE,  we denote
\begin{equation}\label{eq4.5a}
\left\{
\begin{array}{lll} \alpha(t)
 &=&  - b_x(t, X(t), Y(t), \xi(t)   )   -b_y(t, X(t), Y(t), \xi(t) )\nabla_{X(t)} F\\ \\
\be(t)  &=&   - \si_x
 (t, X(t), Y(t),   \xi(t))   -\si_y
 (t, X(t), Y(t),   \xi(t))   \nabla_{X(t)} F\\ \\
\phi(t)& =&    -  f_x(t, X(t), Y(t),   \xi(t)) -f_y(t, X(t), Y(t),
 \xi(t))\nabla_{X(t)} F\\ \\
 \Theta &=&  g_x(X(T), Y(T))+\EE \left[ g_y(X(T), Y(T)) \right] \nabla_{X(T)} F\,.
\end{array}\right.
\end{equation}
We also denote
\[
\rho _{t,T}:= \exp\left\{\int_t^T  \beta(s)  dW_{s}+\int_t^T \left[
\al(s) ds-\frac12 \beta^2(s)  \right] ds\right\}\,.
\]
Then application of   Equation (2.11) of \cite{hns} to  the above
backward SDE yields
\begin{equation}
p(t)  =  \mathbb{E}\left( \Theta \rho _{t,T}+\int_{t}^{T}\rho _{t,r}\phi(r) dr\big|\mathcal{F%
}_{t}\right) \,.  \label{e.sol-lin-bsde-1}
\end{equation}%

%

Substituting this into the maximum principle equations \eqref{eq3.18a}, we see that
the maximum principle consists of the following equations:
\begin{equation}
\frac{\partial f}{\partial \xi} (t, X(t), Y(t),   \xi(t)) +\la(t,
X(t), Y(t)) p(t) + h(t, X(t), Y(t))\le 0 \label{e.domain}
\end{equation}
and
\begin{equation}
\left[\frac{\partial \phi }{\partial \xi} (t, X(t), Y(t),   \xi(t))
+\la(t, X(t), Y(t)) p(t) + h(t, X(t),
Y(t))\right]d\xi(t)=0\,.\label{e.control}
\end{equation}

Thus here is the strategy to solve the maximum principle equations
\eref{e.system}, \eref{e.domain} and \eref{e.control}, where the
$p(t)$ in \eref{e.domain} and \eref{e.control} is given by
\eref{e.sol-lin-bsde-1}: First, assuming that $X$ is known, we solve
the backward stochastic differential equation \eref{e.bsde} with the
terminal condition \eref{e.bsde-t} to obtain
\eref{e.sol-lin-bsde-1}. Then substitute it into \eref{e.domain}.
This way  we obtain an equation which
  will describe the domain $D$ in which
the process $X(t)$ must be all the time. We control the process $X(t)$
 in such a way that when the process is in the
interior of the domain $D$,  we don't do anything. When the process
reaches  the boundary of $D$, we exercise the minimal push to keep the
process inside the domain $D$.  Here are some detailed explanation
of the above strategy:

The equations \eref{e.domain} and \eref{e.control} are essentially
an equation for the ``domain" of the state $X(t)$ and the condition
for the singular control $\xi$ to satisfy. The equation
\eref{e.domain} can be complicated since the solution $p(t)$ may
depend on the paths  of $X$, $Y$ and the path of control $\xi$ itself
up to time $t$.  Denote
$X_t=(X(s), 0\le s\le t)$ the trajectory of $X$ up to time $t$. Then
$p(t)$ can be represented in general as $p(t)=p(X(t), Y(t), \xi(t))$.

We now consider a slightly more general situation,  where the singular
control may be any finite variation process, not necessarily increasing.
The increasing case corresponds to   $r(t, X_t, Y_t, \xi_t)=\infty$ below.

 Suppose that there are two functionals $l, r: [0, T]\times C([0,
 T],\RB)^3\rightarrow \RB$ with $l\le r$
  such that the
equations \eref{e.domain} and \eref{e.control} can be written as
\begin{equation}
\left\{\begin{array}{ll} &l(t, X_t, Y_t, \xi_t) \le X(t)\le r(t,
X_t, Y_t, \xi_t)\,\\  \\
&\int_0^T \left[ X(t)-l(t, X_t, Y_t,
\xi_t)\right]
d\xi(t)=0\\
\\
&\int_0^T \left[r(t, X_t, Y_t, \xi_t)-X(t)\right] d\xi(t)=0\,.
\end{array}\right.\label{e.domain-1}
\end{equation}
Then we are led to   the  problem of  finding  a finite variation
(not necessarily increasing)  control
$\xi$ for the system
\begin{eqnarray}
dX(t) &=& b(t, X(t), Y(t), \xi(t)  ) dt+ \si (t, X(t), Y(t), \xi(t)
) dB(t)\nonumber\\
&&\qquad + \la(t, X(t), Y(t))d\xi(t) \label{e.system-1}
\end{eqnarray}
satisfying \eref{e.domain-1}.  This is a   Skorohod type problem.
For simplicity, we restrict ourselves to the case when
\[
\la(t,x,y)=1\,.
\]

\begin{theorem} Suppose that the following hold
\begin{enumerate}
\item \ $b$ and $\si$ are uniformly Lipschitz continuous.  Namely,
there is a positive  constant $L$ such that
\begin{equation}
|b(t,x_2,y_2,\xi_2)-b(t,x_1,y_1,\xi_1)|\le
L(|x_2-x_1|+|y_2-y_1|+|\xi_2-\xi_1|)\,.
\end{equation}
The same inequality holds for $\si$.
\item \ $l$ and $r$ are uniformly Lipschitz continuous, i.e.
\begin{align}
& |r(t, X_t^2, Y_t^2, \xi_t^2)-r(t, X_t^1, Y_t^1, \xi_t^1)| 
 \le \kappa \sup_{0\le s \le t} \left[ |X^2(s)-X^1(s) |+|\xi^2(s)-\xi^1(s)| \right]
\nonumber\\
& \qquad + L
\int_0^t  \sup_{0\le s \le r} \big[ |X^2(s)-X^1(s)
+|Y^2(s)-Y^1(s)|+|\xi^2(s)-\xi^1(s)| \big] dr
\end{align}
for some $\kappa <1/4$.
The same inequality holds for $l$.
\item \  For any $t\in [0, T]$, $X$, $Y$ and $\xi$ in $C([0,
T]\,, \RB)$,
\begin{equation}
l(t, X_t, Y_t, \xi_t)<r(t, X_t, Y_t, \xi_t)\,.
\end{equation}
\end{enumerate}
Then,  Equations \eref{e.domain-1}-\eref{e.system-1} have a unique
solution.
\end{theorem}

\dproof  We shall apply the Banach fixed point theorem to prove the theorem.
Let us denote by  $\BB$ the Banach space of all
continuous adapted processes $(X(t), \xi(t))$ which are square
integrable. More precisely,
\begin{eqnarray*}
\BB&=&\bigg\{ (X, \xi)\,, \ \hbox{$X$ and $\xi$ are continuous and
adapted \ and }\\
&&\qquad \|(X, \xi)\|_\BB:=  \left\{\EE \sup_{0\le t\le T}\left(
|X(t)|^2+ |\xi(t)|^2\right)\right\}^{1/2}<\infty\bigg\}\,.
\end{eqnarray*}
From \eref{e.domain-1} and \eref{e.system-1},
we define  the following mapping on $\BB$:
\begin{eqnarray}
&&F(X, \xi)=(Z, \eta):\\
&&\hbox{where $(Z,\eta)$ satisfies the inequalities:}\nonumber\\
&&\left\{\begin{array}{ll} &l(t, X_t, Y_t, \xi_t) \le Z(t)\le r(t,
X_t, Y_t, \xi_t)\,\\  \\
&\int_0^T \left[ X(t)-l(t, X_t, Y_t, \xi_t)\right]
d\eta(t)=0\\
\\
&\int_0^T \left[r(t, X_t, Y_t, \xi_t)-X(t)\right] d\eta(t)=0\,,
\end{array}\right.\label{e.domain-2}
\end{eqnarray}

and, in addition,
\begin{eqnarray}
dZ(t) = b(t, X(t), Y(t), \xi(t)  ) dt+ \si (t, X(t), Y(t), \xi(t) )
dB(t)+ d\eta(t) \label{e.system-2}
\end{eqnarray}
For every given continuous pair $(X(t), \xi(t))$ in $\BB$, by the condition (3),
Theorem 2.6   and  Corollary 2.4 of
\cite{bkr} the above Skorohod problem has a unique solution $(Z(t),
\eta(t))$  and the
solution pair $(Z(t), \eta(t))$   can be represented as
\begin{eqnarray}
\eta(t)&=&\Xi(l, r, \psi ) (t)\label{e.eta}\\
Z(t)&=& \psi(t)-\eta(t)\label{e.z}\,,
\end{eqnarray}
where
\begin{eqnarray*}
\psi(t)&:=& \int_0^t b(s, X(s), Y(s), \xi(s)  ) ds+ \int_0^t \si (s,
X(s), Y(s), \xi(s) )
dB(s)\\
l(u)&:=& l(u, X_u, Y_u, \xi_u)\nonumber\\
r(u)&:=& r(u, X_u, Y_u, \xi_u)\nonumber\\
 \Xi(l, r, \psi)  (t) &:=&  \max\Bigg\{ \left[ (\psi(0)-r(0))^+
\wedge
\inf_{u\in [0, t]} (\psi(u)-l(u))\right]\,, \nonumber\\
&&\qquad \sup_{s\in [0, t]} \left[ (\psi(s)-r(s))  \wedge \inf_{u\in
[s, t]} (\psi(u)-l(u))\right]\Bigg\} \,.
\end{eqnarray*}
It is elementary that
\begin{eqnarray*}
&&\left|\max \{b_1, b_2\}-\max \{a_1, a_2\}\right|\le \max \{b_1-a_1, b_2-a_2\}\\
&&\left|\sup _{0\le t\le T}g(t)-\sup _{0\le t\le T}f(t)
\right|\le \sup _{0\le t\le T}|g(t)-f(t)|\\\
&&\left|\inf_{0\le t\le T}g(t)-\inf_{0\le t\le T}f(t)
\right|\le \sup_{0\le t\le T}|g(t)-f(t)|\,.
\end{eqnarray*}
 From the  expression of $\Xi$, we easily see that
\begin{align}
\sup_{0\le s\le t}| \Xi(l_2, r_2, \psi_2)(s)-\Xi(l_1, r_1,
\psi_1)(s)|
\le  & 2 \sup_{0\le s\le t} \left[ |l_2(s)-l_1(s)|+
|r_2(s)-r_1(s)|\right] \nonumber\\
&\qquad + 4 \sup_{0\le s\le t} \left[|\psi_2(s)-\psi_1  (s)|\right]\,.\label{e.xi}
\end{align}
Now we want to show that $\BB \ni (X, \xi) \rightarrow F(X, \xi)=(Z,
\eta)$ is a contraction on $\BB$.  Assume that $(X^1, \xi^1)$ and
$(X^2, \xi^2)$ be two elements in $\BB$ and let $(Z^1, \eta^1)$ and
$(Z^2, \eta^2)$ be the corresponding solutions to
\eref{e.domain-2}-\eref{e.system-2}. Then  for $i=1,2$, we have
\begin{eqnarray}
\eta_i(t)&=&\Xi(l_i, r_i, \psi_i )(t)\label{e.eta-i}\\
Z_i(t)&=& \psi_i(t)-\eta_i(t)\label{e.z-i}\,,
\end{eqnarray}
where
\begin{eqnarray*}
\psi_i(t)&:=& \int_0^t b(s, X^i(s), Y^i(s), \xi^i(s)  ) ds+ \int_0^t
\si_i (s, X^i(s), Y^i(s), \xi^i(s) )
dB(s)\\
l_i(u)&:=& l(u, X_u^i, Y_u^i, \xi_u^i)\nonumber\\
r_i(u)&:=& r(u, X_u^i, Y_u^i, \xi_u^i)\nonumber\\
\eta_i&:=& \Xi(l_i, r_i, \psi_i)(t) =  \max\Bigg\{ \left[ (\psi_i(0)-r_i(0))^+
\wedge
\inf_{u\in [0, t]} (\psi_i(u)-l_i(u))\right]\,, \nonumber\\
&&\qquad \sup_{s\in [0, t]} \left[ (\psi_i(s)-r_i(s))  \wedge
\inf_{u\in [s, t]} (\psi_i(u)-l_i(u))\right]\Bigg\} \,.
\end{eqnarray*}
From \eref{e.xi} and then from the assumptions on $l$ and $r$, we
see that
\begin{eqnarray*}
&&\EE \sup_{0\le r\le t}|\eta_2(r)-\eta_1(r)|^2 \\
&\le&  8 \EE
\sup_{0\le s\le t} \left[ |l_2(s)-l_1(s)|^2+
|r_2(s)-r_1(s)|^2\right] + 32\EE
\sup_{0\le s\le t} \left[ |\psi_2(s)-\psi_1 (s)|^2\right]\\
&\le& 8\kappa^2\EE \sup_{0\le r\le t} \left[
|X^2(r)-X^1(r)|^2+|\xi^2(r)-\xi^1(r)|^2 \right]\\
&&\qquad  +32 \EE \sup_{0\le s\le t} |\psi_2(s)-\psi_1 (s)|^2\\
&&\qquad   +C \int_0^t \EE \sup_{0\le r\le t}\left[
|X^2(r)-X^1(r)|^2+|\xi^2(r)-\xi^1(r)|^2  \right] dr \,.
\end{eqnarray*}
By standard argument from stochastic analysis, we have
\[
\EE \sup_{0\le s\le t} |\psi_2(s)-\psi_1 (s)|^2 \le C \int_0^t \EE
\sup_{0\le s\le r}\left[ |X^2(s)-X^1(s)|^2+|\xi^2(s)-\xi^1(s)|^2
\right] dr\,.
\]
Thus we have
\begin{eqnarray}
 \EE \sup_{0\le r\le t}|\eta_2(r)-\eta_1(r)|^2
 &\le&8 \EE
\sup_{0\le s\le t} \left[ |l_2(s)-l_1(s)|^2+
|r_2(s)-r_1(s)|^2\right] \nonumber\\
&&\qquad +
 C \int_0^t \EE \sup_{0\le s\le r}\left[
|X^2(s)-X^1(s)|^2+|\xi^2(s)-\xi^1(s)|^2  \right] dr
\,.\nonumber\\
\label{e.eta-est}
\end{eqnarray}
From \eref{e.z-i} we have
\begin{eqnarray*}
 \EE \sup_{0\le r\le t}|Z_2(r)-Z_1(r)|^2
 &\le& 2 \EE \sup_{0\le r\le t}|\eta_2(r)-\eta_1(r)|^2+2 \EE \sup_{0\le r\le t}|\psi_2(r)-\psi_1(r)|^2
 \nonumber\\
 &\le&16\kappa^2\EE \sup_{0\le r\le t} \left[
|X^2(r)-X^1(r)|^2+|\xi^2(r)-\xi^1(r)|^2 \right]\\
&&\qquad  +64 \EE \sup_{0\le s\le t} |\psi_2(s)-\psi_1 (s)|^2\\
&&\qquad   +C \int_0^t \EE \sup_{0\le s\le r}\left[
|X^2(s)-X^1(s)|^2+|\xi^2(s)-\xi^1(s)|^2  \right] dr\\
 &\le&16\kappa^2\EE \sup_{0\le r\le t} \left[
|X^2(r)-X^1(r)|^2+|\xi^2(r)-\xi^1(r)|^2 \right]\\
&&\qquad   +C \int_0^t \EE \sup_{0\le s\le r}\left[
|X^2(s)-X^1(s)|^2+|\xi^2(s)-\xi^1(s)|^2  \right] dr
\,.
\end{eqnarray*}
Combining the above inequality with \eref{e.eta-est}, we have
\begin{eqnarray}
 &&\EE \sup_{0\le r\le t}\left [|Z_2(r)-Z_1(r)|^2+|\eta_2(r)-\eta_1(r)|^2\right] \nonumber\\
 &\le& 16\kappa^2\EE \sup_{0\le r\le t} \left[
|X^2(r)-X^1(r)|^2+|\xi^2(r)-\xi^1(r)|^2 \right] \nonumber\\
&&\qquad   +C \int_0^t \EE \sup_{0\le s\le r}\left[
|X^2(s)-X^1(s)|^2+|\xi^2(s)-\xi^1(s)|^2  \right] dr\nonumber\\
&\le& (16\kappa^2+Ct)\EE \sup_{0\le r\le t} \left[
|X^2(r)-X^1(r)|^2+|\xi^2(r)-\xi^1(r)|^2 \right]\,.  
 \label{e.z-est}
\end{eqnarray}
If $\kappa<1/4$,  then we can choose $t_0$ such that  $16\kappa^2+Ct<1$ for all
$t\le t_0$.
Thus from \eref{e.z-est}, we conclude that $F$ is a
contraction mapping from $\BB$ to $\BB$.   Following a routine argument, we see that
the solution $\xi(t), X(t)$ for the equations \eref{e.domain-1} and \eref{e.system-1}
 up to time $t_0$.  Since the constant $C$ in
$16\kappa^2+Ct$  does not depends on the initial condition, we repeat this procedure
to solve the equations  \eref{e.domain-1} and \eref{e.system-1} for on the interval $[0, T]$.
\fproof 
\begin{remark}From the proof of the theorem, we see that if we define the Picard iteration
for $n=0, 1, 2, \cdots$,
\begin{equation}
\left\{\begin{array}{ll} &l(t, X_t^{(n)} , Y_t^{(n)}, \xi_t^{(n)}) \le X^{(n+1)}(t)\le r(t,
X_t^{(n)}, Y_t^{(n)}, \xi_t^{(n)})\,\\  \\
&\int_0^T \left[ X^{(n)}(t)-l(t, X_t^{(n)}, Y_t^{(n)}, \xi_t^{(n)})\right]
d\xi^{(n+1)}(t)=0\\
\\
&\int_0^T \left[r(t, X_t^{(n)}, Y_t^{(n)}, \xi_t^{(n)})-X^{(n)}(t)\right] d\xi^{(n+1)}(t)=0\,,
\end{array}\right.\label{e.domain-3}
\end{equation}
and
\begin{eqnarray}
dX^{(n+1)}(t)
&=& b(t, X^{(n)}(t), Y^{(n)}(t), \xi^{(n)}(t)  ) dt\label{e.system-3}\\
&&\qquad  + \si (t, X^{(n)}(t), Y^{(n)}(t), \xi^{(n)}(t) )
dB(t)+ d\xi^{(n+1)}(t) \,,\nonumber
\end{eqnarray}
where $X_t^{(0)}=X(0)$, $\xi_t^{(0)}=0$, then  $(X^{(n)}(t), \xi^{(n)}(t))$ will converge
to the true solution $(X(t), \xi(t))$ in $\BB$.  This may be used to construct the
numerical solutions.
\end{remark}

\begin{example}\rm 
Let us consider an optimal harvesting problem where the population density $X(t)$ at time $t$ is described by the linear controlled system
\begin{eqnarray}\label{ex4.3}
dX(t) &=& \left[ b_0(t, \xi)+b_1(t) X(t)+b_2(t) \EE X(t)\right]
dt\nonumber\\
&&\qquad  + \si(t, X(t), Y(t), \xi(t)) dB(t) -  d\xi(t)\,.
\label{e.lin-system}
\end{eqnarray}
We allow the coefficients $b_0$ and $\sigma_0$ to depend on the harvested amount $\xi$ to model the situation where the harvesting has influence on the environment and hence on the population growth.
We want to find $\hat \xi$ such that
\begin{equation}
\sup_{\xi\in \cA} J(\xi)=J(\hat \xi)\,,\label{e.4.29}
\end{equation}
where
\[
J(\xi)=\EE \left[\int_0^T f(t, X(t), Y(t), \xi(t)) dt+g(X(T), Y(T))+\int_0^T h(t, X(t))  \xi(dt)\right]\,,
\]
with
\begin{equation}
f(t, x, y,\xi)=f_1(t)x+f_2(t)y+f_3(t, \xi)
\label{e.4.27}
\end{equation}
and
\begin{equation}
g(x,y)=Kx\label{e.4.28}
\end{equation}
with $K>0$.
Then from \eqref{eq4.5a} we get 
\[
\left\{
\begin{array}{ll} \al(t)
&=  - b_1(t)-b_2(t) \\ \\
\be(t)  &=   - \si_x
 (t, X(t), Y(t),   \xi(t))   -\si_y
 (t, X(t), Y(t),   \xi(t))   \nabla_{X(t)} F\\ \\
\phi(t)&=    -  f_1(t)-f_2(t) \\ \\
 \Theta &=  K\,.
 \end{array}\right.
\]
Denote
\[
\rho _{t,T}= \exp\left\{\int_t^T  \beta(s)  dW_{s}+\int_t^T \left[
\al(s) ds-\frac12 \beta^2(s)  \right] ds\right\}\,.
\]
 Since $\al$ is deterministic,  we have
for all $t\le r\le T$,
\begin{eqnarray*}
\al(t,r)&:=&\mathbb{E}\left(   \rho _{t,r} \big|\mathcal{F%
}_{t}\right)\\
 &=& \exp\left\{\int_t^r \al(s) ds\right\} \mathbb{E}\left(   \exp\left\{\int_t^r  \beta(s)
dW_{s}-\frac12 \int_t^r    \beta^2(s)
ds\right\} \big|\mathcal{F%
}_{t}\right)
\\
&=&\exp\left\{\int_t^r \al(s) ds\right\}\,.
\end{eqnarray*}
Note that $\al(t,r)$ is a deterministic function. It is easy  to see from
\eref{e.sol-lin-bsde-1}  that
\begin{equation}
p(t)  =  K \al({t,T})+\int_{t}^{T}\al({t,r}) \phi (r) dr\label{e.4.30}
\end{equation}%
is a deterministic function.
Thus we have
\[
\frac{\partial}{\partial \xi}f_3(t,\xi) +\la(t, X(t), Y(t))p(t)+h(t, X(t), Y(t))\le 0\,.
\]
\[
\left[\frac{\partial}{\partial \xi}f_3(t,\xi) +\la(t, X(t), Y(t))p(t)+h(t, X(t),
Y(t))\right]d\xi(t)= 0\,.
\]
If furthermore we assume
\begin{equation}
\phi(t)\ge 0\,, \quad \frac{\partial}{\partial \xi}f_3(t,\xi)=0\,, \quad
{\rm and}\quad h(t,x,y)=h_0(t) x^{\kappa }\,, \label{e.4.31}
\end{equation}
where  $h_0(t)$ is positive and $\kappa$ is a constant,  we get  (noting that $\la(t)=-1$)
\[
-  p(t) +h_0(t)X^{\kappa }(t)\le 0\,,
\]
or
\[
X(t)\begin{cases}
 \le  \left(\frac{p(t)}{h_0(t)}\right)^{\frac1{\kappa }}& \qquad \hbox{if}\ \kappa>0\\ \\
 \ge  \left(\frac{ h_0(t)}{p(t)}  \right)^{-\frac1{\kappa}}& \qquad \hbox{if}\ \kappa<0 \,. \\
\end{cases}
\]
In this case, we can take
\begin{equation}
\left\{
\begin{array}{lll}
l(t,x,y,\xi)=0\qquad {\rm and}\qquad
r(t,x,y,\xi) = \left(\frac{p(t)}{h_0(t)}\right)^{\frac1{\kappa }}& \qquad \hbox{if}\ \kappa>0\\ \\
l(t,x,y,\xi)=\left(\frac{ h_0(t)}{p(t)}  \right)^{-\frac1{\kappa}}\qquad {\rm
and}\qquad
r(t,x,y,\xi)=\infty & \qquad \hbox{if}\ \kappa<0 \,. \\
\end{array}\right. \label{e.4.32}
\end{equation}
Note that $\kappa<0$ means that unit price goes up when the population goes down
(which becomes more precious).
In this case, we want keep the population above a threshold
$\bar h(t)=\left(\frac{ h_0(t)}{p(t)}  \right)^{-\frac1{\kappa}}$. It is interesting to note
that when  $h_0(t)$ is larger, this threshold $\bar h_0(t)$ is also larger.
We have proved
\begin{theorem}\label{t.4.4}
Under the assumptions \eref{e.4.27}, \eref{e.4.28} and \eref{e.4.31}, the solution $\hat \xi$
of the mean field singular control  problem \eref{e.4.29} is given by the solution $(\hat X,
\hat Y, \hat \xi)$ of the Skorohod reflection problem \eref{e.domain-1} and \eref{e.lin-system},
with the boundaries $l$ and $r$ given by \eref{e.4.32}.
\end{theorem}

Next,  we continue the above example but with $h$ being given by
\begin{equation}
h(t,x,y):=h_0(t)x^2+h_1(t) x\,.\label{e.h-quadratic}
\end{equation}
Namely, we continue to assume \eref{e.lin-system}-\eref{e.4.30}.  But we replace \eref{e.4.31} by
\begin{equation}
\phi(t)\ge 0\,, \quad \frac{\partial}{\partial \xi}f_3(t,\xi)=0\,, \quad
{\rm and}\quad h(t,x,y)=h_0(t)x^2+h_1(t) x\,. \label{e.4.33}
\end{equation}
where  $h_0(t)$ is positive.  Then,  the inequalities  \eref{e.domain} -  \eqref{e.control} become
\[
-p(t) +h_0(t)X^2(t) +h_1(t) X(t)\le 0\,.
\]
\begin{equation}
l(t)\le X(t)\le r(t)\,,\label{e.4.34}
\end{equation}
where
\begin{eqnarray}
l(t)
&:=& \frac{h_1(t)-\sqrt{h_1^2(t)+4h_0(t) p(t)}} {2h_0(t)}\label{e.4.35}\\
r(t)
&:=& \frac{h_1(t)+\sqrt{h_1^2(t)+4h_0(t) p(t)}} {2h_0(t)}\,. \label{e.4.36}
\end{eqnarray}
Similar to Theorem \ref{t.4.4}, we have
\begin{theorem}\label{t.4.5}
Under the assumptions \eref{e.4.27}, \eref{e.4.28} and \eref{e.4.33}, the solution $\hat \xi$
of the mean field singular control  problem \eref{e.4.29} is given by the solution $(\hat X,
\hat Y, \hat \xi)$ of the Skorohod reflection problem \eref{e.domain-1} and \eref{e.lin-system},
with the boundaries $l$ and $r$ given by \eref{e.4.34}-\eref{e.4.36}.
\end{theorem}
We can also consider the case that $h$ is given by \eref{e.h-quadratic}
but with $h_0(t)<0$.  In this case,  the domain \eref{e.domain} will be either
\[
X(t)\le \underbar h(t) \quad\hbox{or} \quad
 X(t)\ge \bar h(t)\,.
 \]
The interested readers may write down similar result for this case as well.

\end{example}

\section{General singular mean-field  games}\label{sec3}

\subsection{Statement of the problem}
In this section we consider the stochastic game of two players, each of them
is to maximize his/her singular mean-field performance.

Denote  $\xi = (\xi_1, \xi_2), u=(u_1,u_2), w = (w_1, w_2), \lambda=(\lambda_1,\lambda_2),
h=(h_1,h_2)$ with $h_i=(h_{i,1},h_{i,2})$, and
let the pair $w_i= (\xi_i,u_i)$ represent the control of player $i \; ; \; i=1,2$.\\

Suppose the process $X(t) = X^{\xi, u}(t)$ under control of the two players satisfy the following stochastic differential equation with jumps.
\begin{align}\label{eq3.1}
dX(t) & = b(t,X(t),Y(t),\xi(t), u(t),\omega)dt + \sigma(t,X(t),Y(t), \xi(t),  u(t),\omega)dB(t) \nonumber \\
&
 + \lambda(t,X(t),u(t),\omega)d\xi(t),
\end{align}
where
\begin{equation}\label{eq3.2}
Y(t) = F(X(t, \cdot)),
\end{equation}
and $F$ is a Fr\'echet differentiable operator on $L^2(P)$.

We put $\mathbb{G}^{i} = \{ \mathcal{G}_t^{i} \}_{t\geq 0}$ where $\mathcal{G}_t^{i}\subseteq \FC_t$ is the information available to player $i$ at time $t$.The \emph{performance functional} for player $i$ is assumed to be on the form
\begin{align}\label{eq3.3}
J_i(\xi,u) & = E \left[ \int_0^T f_i(t,X(t),Y(t),w(t),\omega)dt + g_i(X(T),Y(T),\omega)\right. \nonumber \\
 & + \left. \int_0^T h_i(t,X(t),u(t),\omega) \xi(dt)\right] ; \; \; i=1,2.
 \end{align}

  We want to find a \emph{Nash equilibrium} for this game, i.e. find $(\xi^*_1, u^*_1)
   \in \AC^{(1)}$ and $(\xi^*_2, u^*_2) \in \AC^{(2)}$ such that
 \begin{equation}\label{eq3.4}
 \sup_{(\xi_1,u_1)  \in \AC^{(1)}} J_1(\xi_1,u_1,\xi^*_2,u^*_2) = J_1(\xi^*_1,u^*_1,\xi^*_2,u^*_2)
 \end{equation}
 and
 \begin{equation}\label{eq3.5}
 \sup_{(\xi_2,u_2) \in \AC^{(2)}} J_2(\xi^*_1,u^*_1,\xi_2,u_2) = J_2(\xi^*_1,u^*_1,\xi^*_2,u^*_2)
 \end{equation}

 Here $\AC^{(i)}$ is a given family of $\GB^{(i)}$-predictable processes such that the corresponding state equation has a unique solution $X$ such that $\omega \rightarrow X(t,\omega) \in L^2(P)$ for all t.
We let $A^{(i)}$ denote the set of possible values of $u_i(t); t \in [0,T] $ when $(\xi_i,u_i) \in \AC^{(i)}; i= 1,2.$

%

 \subsection{A sufficient maximum principle for the general non-zero sum case}

 Define two {\it Hamiltonians} $H_i; i=1,2,$ as follows:
 \begin{align}\label{eq3.6}
 H_i& (t,x,y,\xi_1,u_1,\xi_2,u_2,p_i,q_i)(dt, \xi_1(dt), \xi_2(dt)) \nonumber \\
 & = H_{i,0} (t,x,y,\xi_1,u_1,\xi_2,u_2,p_i,q_i) dt + \sum_{j=1}^2 \{ \lambda_j (t,x,u) p_i + h_{i,j}(t,x,u)\} \xi_j(dt)
 \end{align}
 where
 \begin{align}\label{eq3.7}
 H_{i,0}& (t,x,y,w,p_i,q_i):= f_i(t,x,y,w) + b(t,x,y,\xi, u)p_i + \sigma(t,x,y,\xi, u)q_i.
 \end{align}
 We assume that for $i=1,2$, $H=H_i$ is Fr\'echet differentiable $(C^1)$ in the variables $x,y,\xi,u$.


 The BSDE for the adjoint processes $p_i,q_i$ is
 \begin{equation}\label{eq3.8}
 \begin{cases}
 dp_i(t) & = \displaystyle - \frac{\partial H_{i,0}}{\partial x} (t,X(t),Y(t),w(t),p_i(t),q_i(t))dt\\
 &\displaystyle- \frac{\partial H_{i,0}}{\partial y} (t,X(t), Y(t), w(t),p_i(t),q_i(t))\nabla_{X(t)} F)dt\\
 &\displaystyle+ q_i(t) dB(t)\\
 p_i(T) & = \displaystyle \frac{\partial g_i}{\partial x}(X(T),Y(T))+\frac{\partial g_i}{\partial y}(X(T),Y(T))\nabla_{X(T)} F; \; \; i= 1,2.\\
 \end{cases}
 \end{equation}
 Note that \eqref{eq3.8} is an operator valued BSDE for each $\omega$.\\

 \begin{theorem}[Sufficient maximum principle]\label{th3.1}
 Let $(\hxi_1,\hu_1) \in \AC^{(1)},(\hxi_2,\hu_2) \in \AC^{(2)}$ with corresponding solutions $\hX,\hp_i,\hq_i,\hr_i$ of \eqref{eq2.1} and \eqref{eq2.6}. Assume the following:
 \begin{itemize}

 \item
The maps
 \begin{equation}\label{eq3.9}
 X,w_1 \rightarrow H_1(t,X,F(X),w_1,\hw_2(t),\hp_1(t),\hq_1(t)),
 \end{equation}
 and
 \begin{equation}\label{eq3.10}
 X,w_2 \rightarrow H_2(t,X,F(X),\hw_1,w_2(t),\hp_2(t),\hq_2(t)),
 \end{equation}
 and
 \begin{equation}\label{eq3.10a}
 X \rightarrow g_i(X,F(X))
 \end{equation}
 are concave for all $t; i=1,2$.

\item (The conditional maximum properties)
 \begin{align}\label{eq3.11}
\ess_{u_1\in A_1} & E [H_{1}(t,\hX(t),\hY(t),\hxi_1(t),u_1, \hxi_2(t),\hu_2(t), \hp_1(t),\hq_1(t)) \mid \GC^{(1)}_t] \nonumber \\
& = E[H_{1}(t,\hX(t), \hY(t),\hxi_1(t),\hu_1(t),\hxi_2(t),\hu_2(t), \hp_1(t),\hq_1(t)) \mid \GC^{(1)}_t]
\end{align}
and
\begin{align}\label{eq3.12}
\ess_{u_2\in A_2} & E[H_{2}(t,\hX(t), \hY(t),\hxi_1(t),\hu_1(t),\hxi_2(t),u_2, \hp_2(t), \hq_2(t)) \mid \GC_t^{(2)}] \nonumber \\
& = E[H_{2}(t,\hX(t), \hY(t),\hxi_1(t),\hu_1(t),\hxi_2(t),\hu_2(t), \hp_2(t), \hq_2(t)) \mid \GC_t^{(2)}],
\end{align}

\item (Variational inequalities)
 \begin{align}\label{eq3.13a}
\ess_{\xi_1} & E [H_{1}(t,\hX(t),\hY(t),\xi_1,\hu_1(t), \hxi_2(t),\hu_2(t), \hp_1(t),\hq_1(t)) \mid \GC^{(1)}_t] \nonumber \\
& = E[H_{1}(t,\hX(t), \hY(t),\hxi_1(t),\hu_1(t),\hxi_2(t),\hu_2(t), \hp_1(t),\hq_1(t)) \mid \GC^{(1)}_t]
\end{align}
and
\begin{align}\label{eq3.13b}
\ess_{\xi_2} & E[H_{2}(t,\hX(t), \hY(t),\hxi_1(t),\hu_1(t),\xi_2,\hu_2(t), \hp_2(t), \hq_2(t)) \mid \GC_t^{(2)}] \nonumber \\
& = E[H_{2}(t,\hX(t), \hY(t),\hxi_1(t),\hu_1(t),\hxi_2(t),\hu_2(t), \hp_2(t), \hq_2(t)) \mid \GC_t^{(2)}],
\end{align}

 \end{itemize}

  Then $(\hxi_1, \hu_1), (\hxi_2, \hu_2)$ is a Nash equilibrium, in the sense that \eqref{eq3.4} and \eqref{eq3.5} hold with $\xi^*_i := \hxi_i, u^*_i := \hu_i ;$ $i=1,2$.
 \vskip 0.3cm

 \end{theorem}

 \dproof 
 By introducing a suitable increasing sequence of stopping times converging to $T$,
  we see that we may assume that all local martingales appearing in the proof below are martingales.
  We refer to \cite{OS6} for details.
 We first study the stochastic control problem \eqref{eq3.4}.
 For simplicity of notation, in the following we put $X(t) = X^{\xi_1,u_1,\hxi_2,\hu_2}(t)$, $Y(t) = Y^{\xi_1,u_1,\hxi_2,\hu_2}(t)$ and $\hX(t) = X^{\hxi_1,\hu_1,\hxi_2,\hu_2}(t)$, $\hY(t) = Y^{\hxi_1,\hu_1,\hxi_2,\hu_2}(t), b(t)= b(t,X(t),Y(t),\xi_1(t),\xi_2(t),u_1(t),u_2(t),\omega), \hb(t)= b(t,\hX(t),\hY(t),\xi_1(t),\hxi_2(t),u_1(t),\hu_2(t),\omega)$ and similarly with $\sigma(t),\hat{\sigma}(t)$.

 Consider $J_1(\xi_1,u_1,\hxi_2,\hu_2) - J_1(\hxi_1, \hu_1,\hxi_2,\hu_2) = I_1 + I_2 + I_3 + I_4$, where
 $$I_1 := E \left[ \int_0^T \{f_1(t,X(t),Y(t),w_1(t),\hw_2(t)) - f_1(t,\hX(t),\hY(t), \hw_1(t), \hw_2(t))\} dt \right]$$
 $$I_2  := E[g_1(X(T),Y(T)) - g_1(\hX(T),\hY(T))]$$
 $$I_3 := E \left[ \int_0^T \{ h_1(t,X(t), u_1(t),\hu_2(t)) d\xi_1(t) - h_{1,1}(t,\hX(t),\hu(t)) d\hxi_1(t)\} \right]$$
$$I_4 := E \left[ \int_0^T \{ h_2(t,X(t),u_1(t),\hu_2(t)) d\xi_2(t) - h_{1,2}(t,\hX(t),\hu(t)) d\hxi_2(t)\} \right]$$

By the definition of $H_1$ we have
\begin{align}\label{eq3.14}
I_1 & = E [ \int_0^T \{ H_{1,0}(t,X_t,Y_t,w_1(t),\hw_2(t),\hp_1(t),\hq_1(t)) \nonumber\\
&\qquad
- H_{1,0}(t,\hX_t,\hY_t,\hw(t), \hp_1(t), \hq_1(t))  - (b - \hb)\hp_1 - (\sigma - \hsi)\hq_1\} dt ]
\end{align}

By concavity of $g_1$ and the It\^o formula we have
\begin{align}\label{eq3.15}
I_2 & \leq E\left[\frac{\partial g_1}{\partial x}(\hX(T),\hY(T)))
 (X(T) - \hX(T))\right.\nonumber\\
 &\qquad\quad  \left. +\frac{\partial g_1}{\partial y}(\hX(T),\hY(T)) \langle\nabla_{\hX(T)} F,X(T) - \hX(T)\rangle\right] \nonumber\\
& =<\hp_1(T),X(T) - \hX(T)> \nonumber\\
& = E \left[ \int_0^T \hp_1(t) d \tX(t) + \int_0^T \tX(t) d\hp_1(t) + \int_0^T \hq_1(t) \tsi(t)dt \right]\,,
\end{align}
where we have put
\begin{equation}\label{eq3.16}
\tX(t) := X(t) - \hX(t), \; \; \tsi(t) := \sigma(t) - \hsi(t). \\
\end{equation}

Note that
\begin{align}\label{eq3.17}
E&\left[ \int_0^T \hp_1(t) d\tX(t)\right] =
 E \left[ \int_0^T \hp_1(t)(b-\hb)dt \right. \nonumber\\
& \left.+ \int_0^T \hp_1(t) (\lambda_1 d\xi_1(t) - \hla_1 d\hxi_1(t) + \lambda_2 d\hxi_2(t) - \hla_2 d\hxi_2(t))\right].
 \end{align}
 and that
 \begin{align}\label{eq3.18}
 E&\left[ \int_0^T \tX(t) d\hp_1(t)\right] =
 E\left[ \int_0^T \tX(t) \{ - \frac{\partial H_{1,0}}{\partial x} (t,X(t), Y(t),w(t),\hp_1(t),\hq_1(t),\hr_1(t)) \right.\nonumber\\
 & \left. - \frac{\partial H_{1,0}}{\partial y} (t,X(t), Y(t), w(t),\hp_1(t),\hq_1(t),\hr_1(t))\nabla_{X(t)} F\}dt \right] 
  \end{align}
Combining \eqref{eq3.18} with $I_3$ and $I_4$ we get
\begin{align}\label{eq3.19}
& \left. J_1(\xi_1,u_1, \hxi_2,\hu_2) - J_1(\hxi_1, \hu_1, \hxi_2, \hu_2) \right.\nonumber\\
& \left. \leq E \left[\int_0^T \{H_1(t) -\hH_1(t) -\frac{\partial \hH_1}{\partial x}(X-\hX) -\frac{\partial \hH_1}{\partial y}\nabla F (X - \hX)\}dt \right]\right. \nonumber\\
& \left. = E \left[ \int_0^T E[H_1(t) -\hH_1(t) -\frac{\partial \hH_1}{\partial x}(X-\hX) -\frac{\partial \hH_1}{\partial y}\nabla F (X - \hX)\mid \GC^{(1)}_t]dt \right]\right. 
\end{align}
where $\hH_1(t)$ means that $H_1$ is evaluated at $(t,\hX(t), \hY(t),\hxi(t),\hu(t),\hp_1(t),\hq_1(t))$, while $H_1(t)$ means that $H_1$ is evaluated at $(t,X(t),Y(t), \xi_1(t), u_1(t),\hxi_2(t),\hu_2(t),\hp_1(t),\hq_1(t))$.

Note that by concavity of $H_1$ we have
\begin{align}\label{eq3.20}
& \left. H_1(t,X,F(X),\xi_1,u_1,\hxi_2,\hu_2,\hp_1,\hq_1) - H_1(t,\hX,F(\hX),\hxi_1,\hu_1,\hxi_2,\hu_2,\hp_1,\hq_1) \right.\nonumber\\
& \left. \leq \frac{\partial \hH_1}{\partial x}(\hX)(X - \hX) +\frac{\partial \hH_1}{\partial y}(\hX) \nabla_{\hX}F(X-\hX)
+ \nabla_{\xi_1} \hH_1(\hxi) (\xi_1 -\hxi_1) + \frac {\partial \hH_1}{\partial u_1}(\hu)(u_1 -\hu_1) \right.
\end{align}
Therefore, to obtain that $J_1 -\hat{J}_1 ≤ 0$,
 it suffices that
\begin{align}\label{eq3.21}
E[\nabla_{\xi_1} \hH_1(\hxi)\mid \GC^{(1)}_t] (\xi_1 -\hxi_1)  \leq 0
\end{align}
for all $\xi_1$, and that\\
\begin{align}\label{eq3.22}
E[\frac {\partial \hH_1}{\partial u_1}(\hu)\mid \GC^{(1)}_t](u_1 -\hu_1) \leq 0
\end{align}
for all $u_1$.
The inequality  \eqref{eq3.22} holds by our assumption \eqref{eq3.11}, and the inequality \eqref{eq3.21} holds by our assumption \eqref{eq3.13a}.
The difference
$$J_2(\hxi_1,\hu_1,\xi_2,u_2) - J_2(\hxi_1,\hu_1, \hxi_2,\hu_2)$$
is handled similarly.
\fproof

 \subsection{The zero-sum game case}

In the {\it zero-sum case} we have
\begin{equation}\label{eq4.1}
J_1(w_1,w_2) + J_2(w_1,w_2) = 0.
\end{equation}
Then the Nash equilibrium $(\hw_1, \hw_2) \in \AC_1 \times \AC_2$ satisfying \eqref{eq3.4}-\eqref{eq3.5} becomes a {\em saddle point} for
\begin{equation}\label{eq4.2}
J(w_1,w_2) := J_1(w_1,w_2).
\end{equation}
To see this, note that \eqref{eq3.4}-\eqref{eq3.5} imply that
$$  J_1(w_1, \hw_2) \leq J_1(\hw_1, \hw_2) = - J_2(\hw_1, \hw_2) \leq -J_2(\hw_1, w_2)$$ and hence
$$  J(w_1, \hw_2) \leq J(\hw_1, \hw_2)  \leq J(\hw_1, w_2) \text{ for all } w_1, w_2.$$
 From this we deduce that
\begin{align}\label{eqA29}
\inf_{w_2 \in \AC_2} &\sup_{w_1 \in \AC_1} J(w_1, w_2) \leq \sup_{w_1 \in \AC_1} J(w_1, \hw_2) \leq J(\hw_1, \hw_2) \nonumber \\
& \leq  \inf_{w_2 \in \AC_2} J(\hw_1, w_2) \leq \sup_{w_1 \in \AC_1} \inf_{w_2 \in \AC_2} J(w_1,w_2).
\end{align}
Since we always have $\inf \sup \geq \sup \inf$, we conclude that
\begin{align}\label{eqA29eg}
\inf_{w_2 \in \AC_2} &\sup_{w_1 \in \AC_1} J(w_1, w_2) = \sup_{w_1 \in \AC_1} J(w_1, \hw_2) = J(\hw_1, \hw_2) \nonumber \\
& =  \inf_{w_2 \in \AC_2} J(\hw_1, w_2) = \sup_{w_1 \in \AC_1} \inf_{w_2 \in \AC_2} J(w_1,w_2).
\end{align}
i.e. $(\hw_1, \hw_2) \in \AC_1 \times \AC_2$ is a {\em saddle point} for  $J(w_1,w_2) $.

Hence we want to find $(\xi^*, \theta^*) \in \AC_\EC \times \Theta$ such that
\begin{equation}\label{eq4.5}
\sup_{\xi \in \AC_\EC} \left\{ \inf_{\theta \in \Theta} J(\xi, \theta)\right\} = \inf_{\theta \in \Theta} \left\{ \sup_{\xi \in \AC_\EC} J(\xi, \theta)\right\} = J (\xi^*, \theta^*),
\end{equation}
where
\begin{align}\label{eq4.6}
J(\xi,u) & = E \left[ \int_0^T f(t,X(t),Y(t),w(t),\omega)dt + g(X(T),Y(T),\omega)\right. \nonumber \\
 & + \left. \int_0^T h(t,X(t),u(t),\omega) \xi(dt)\right].
 \end{align}
As shown in \cite{OS4}, in this case only one Hamiltonian H is needed, namely
\begin{align}\label{eq4.7}
 H& (t,x,y,\xi_1,u_1,\xi_2,u_2,p,q)(dt, \xi_1(dt), \xi_2(dt)) \nonumber \\
 & = H_{0} (t,x,y,\xi_1,u_1,\xi_2,u_2,p,q) dt + \sum_{j=1}^2 \{ \lambda_j(t,x,u) p + h_j(t,x,u)\} \xi_j(dt)
 \end{align}
 where
 \begin{align}\label{eq4.8}
 H_0:=H_{1,0}& (t,x,y,w,p,q)\nonumber \\
 & = f(t,x,y,w) + b(t,x,y,u)p + \sigma(t,x,y,u)q
 \end{align}
and we have put $g_i=g, h_i = h_{1,i} \; ; \; i=1,2$ and $f_1 = f = - f_2$.

Moreover, there is only one couple $(p,q)$ of adjoint processes, given by the BSDE
 \begin{equation}\label{eq4.9}
 \begin{cases}
 dp(t) & = \displaystyle - \frac{\partial H_{0}}{\partial x} (t,X(t),Y(t),w(t),p(t),q(t))dt\\
 &\displaystyle- \frac{\partial H_{0}}{\partial y} (t,X(t), Y(t), w(t),p(t),q(t))\nabla_{X(t)} F)dt\\
  &\displaystyle+ q(t) dB(t)\\

 p(T) & = \displaystyle \frac{\partial g}{\partial x}(X(T),Y(T))+E[\frac{\partial g}{\partial y}(X(T),Y(T))]\nabla_{X(T)} F.\\
 \end{cases}
 \end{equation}

We can now state the corresponding sufficient maximum principle for the
zero-sum game:

\begin{theorem}[Sufficient maximum principle for zero-sum singular mean-field games]\label{th4.1}
Let $(\hw_1, \hw_2) \in \AC_1 \times \AC_2$, with corresponding solutions $\hX(t), \hY(t), \hp(t),\hq(t)$.
Suppose the following holds
\begin{itemize}
\item
The function
 \begin{equation}
 X,w_1 \rightarrow H(t,X,F(X),w_1,\hw_2(t),\hp(t),\hq(t))
 \end{equation}
 is concave for all t, the function
 \begin{equation}
 X,w_2 \rightarrow H(t,X,F(X),\hw_1(t),w_2,\hp(t),\hq(t))
 \end{equation}
 is convex for all t, and the function
 \begin{equation}
 X \rightarrow g(X,F(X))
 \end{equation}
 is affine.\\

\item (The conditional maximum property)
\begin{align}\label{eq4.12}
\ess_{v_1 \in A_1} E [H(t, \hX(t), \hY(t),\hxi_1(t),v_1, \hxi_2(t),\hu_2(t), \hp(t), \hq(t)) \mid \GC_t^{(1)}] \nonumber \\
= E[H(t, \hX(t), \hY(t), \hxi_1(t), \hu_1(t), \hxi_2(t),\hu_2(t), \hp(t), \hq(t))\mid \GC_t^{(1)}]
\end{align}
and
\begin{align}\label{eq4.13}
\esi_{v_2 \in A_2} E [H(t, \hX(t), \hY(t), \hxi_1(t),\hu_1(t), \hxi_2(t),v_2, \hp(t), \hq(t)) \mid \GC_t^{(2)}] \nonumber \\
= E[H(t, \hX(t), \hY(t), \hxi_1(t), \hu_1(t), \hxi_2(t), \hu_2(t), \hp(t), \hq(t))\mid \GC_t^{(2)}].
\end{align}

\item (Variational inequalities)
\begin{align}\label{eq4.14}
\ess_{\xi_1} E [H(t, \hX(t), \hY(t),\xi_1,\hu_1(t), \hxi_2(t),\hu_2(t), \hp(t), \hq(t)) \mid \GC_t^{(1)}] \nonumber \\
= E[H(t, \hX(t), \hY(t), \hxi_1(t), \hu_1(t), \hxi_2(t),\hu_2(t), \hp(t), \hq(t))\mid \GC_t^{(1)}]
\end{align}
and
\begin{align}\label{eq4.15}
\esi_{\xi_2} E [H(t, \hX(t), \hY(t), \hxi_1(t),\hu_1(t), \xi_2,\hu_2(t), \hp(t), \hq(t)) \mid \GC_t^{(2)}] \nonumber \\
= E[H(t, \hX(t), \hY(t), \hxi_1(t), \hu_1(t), \hxi_2(t), \hu_2(t), \hp(t), \hq(t))\mid \GC_t^{(2)}].
\end{align}

\end{itemize}

Then $\hu(t) = (\hu_1(t), \hu_2(t))$ is a saddle point for $J(u_1,u_2)$.
\end{theorem}
\dproof  The proof is similar to (and simpler than) the proof of Theorem \ref{th3.1} and is omitted.
\fproof

\subsection{A necessary maximum principle for the general case}\label{sec5}

In Section \ref{sec6.2} we proved a verification theorem, stating that
if a given control $(\hxi, \hu)$ satisfies certain conditions, then it is
indeed optimal for the singular control game. We now establish a partial converse,
implying that it a control $(\hxi, \hu)$ is optimal for the singular control game,
 then it is a conditional saddle point for the Hamiltonian.

\begin{theorem}\label{th5.1}{(Necessary maximum principle for singular mean-field games)}

Suppose $\hw_1 = (\hxi_1, \hu_i) \in \AC^{(1)}$ and $\hw_2 = (\hxi_2, \hu_2) \in \AC^{(2)}$ constitute a Nash equilibrium for the game, i.e. satisfies \eqref{eq3.4} and \eqref{eq3.5}. Then
\begin{equation}\label{eq5.7}
E \left[ \frac{\partial H_{1,0}}{\partial u_1} (t, \hX(t), \hY(t), \hxi_1(t), u_1, \hxi_2(t), \hu_2(t), \hp_1(t), \hq_1(t))_{u_1=\hu_1(t)} \mid \GC^{(1)}_t \right] = 0
\end{equation}
and
\begin{equation}\label{eq5.8}
E \left[ \frac{\partial H_{2,0}}{\partial u_2} (t, \hX(t), \hY(t), \hxi_1(t), \hu_1(t), \hxi_2(t), u_2, \hp_2(t), \hq_2(t))_{u_2=\hu_2(t)} \mid \GC^{(2)}_t \right] = 0.
\end{equation}
Moreover, the following variational inequalities hold:
\begin{equation}\label{eq5.9}
\begin{cases}
\displaystyle E \left[ \frac{\partial f_i}{\partial \xi_i} (t, \hX(t), \hY(t), \hw(t)) + \lambda_i(t, \hX(t), \hu(t)) \hp_i(t) + h_{ii} (t, \hX(t), \hu(t)) \mid \GC^{(i)}_t \right] \leq 0 \\
\text{ for all } t, \; i=1,2 \\
\text{and} \\
\displaystyle E \left[\frac{\partial f_i}{\partial \xi_i} (t, \hX(t), \hY(t), \hw(t)) + \lambda_i(t, \hX(t), \hu(t)) \hp_i(t) + h_{ii} (t, \hX(t), \hu(t)) \mid \GC^{(i)}_t \right] d \hxi_i(t) =  0 \\
\text{ for all } t, \; i=1,2.
\end{cases}
\end{equation}
\end{theorem}

\dproof  This theorem can be proved in a way similar to the proof  of
Theorem \ref{th5.1a} with an adjustment to the stochastic game case.
The adjustment is similar to one in the proof of Theorem \ref{th3.1}.
\fproof

\section{Applications}\label{sec6}

\subsection{Return to the optimal harvesting problem}\label{sec6.1}

To illustrate our results, we apply it to the optimal harvesting problem \eqref{eq2.1},\eqref{eq2.2} in Section 2.1:
Here the Hamiltonian \eqref{eq5.2a} gets the form

\begin{equation} \label{eq6.1}
H(t,x,y,\xi,p,q)(dt,d\xi) = [yb(t)p + x\sigma(t)q]dt + \{-\lambda_0(t)p + xh_0(t)\} d\xi(t)
\end{equation}

By Theorem 5.1 the corresponding BSDE reduces to
\begin{equation}\label{eq6.2}
 \begin{cases}
 dp(t) & = \displaystyle - [b(t)p(t) + \sigma(t) q(t)]dt + q(t) dB(t)\\
  p(T) & = K\\
\end{cases}
\end{equation}
with associated variational inequalities
\begin{equation}\label{eq6.3}
\begin{cases}
&  -\lambda_0(t)p(t) + h_0(t)X(t) \leq 0;  \; \;  t \in [0,T] \\
 &  [-\lambda_0(t)p(t) + h_0(t)X(t)] d\hxi(t) = 0; \; \;  t \in [0,T]
 \end{cases}
 \end{equation}

 The system consisting of \eqref{eq2.1},\eqref{eq2.2} combined with \eqref{eq6.2},\eqref{eq6.3} represents a \emph{mean-field forward-backward reflected SDE }. 
 Combining the above with the result of Section 4 we get:

 \begin{theorem}\label{th6.1}
Assume that $\hX(t), \hp(t),\hq(t),\hxi(t)$ is a solution of the system \eqref{eq2.1},\eqref{eq2.2} \& \eqref{eq6.2},\eqref{eq6.3}. Then $\hxi(t)$ is an optimal harvesting strategy for the problem \eqref{eq2.3}. Heuristically the optimal harvesting strategy can be described as follows:
 \begin{itemize}
 \item
 If
 \begin{equation}
 \hp(t) > \frac{h_0(t)\hX(t)}{\lambda_0(t)}, \text{i.e.}\; \; \hX(t) < \frac{\lambda_0(t) \hp(t)}{h_0(t)},
 \end{equation}
 then do nothing (choose $\hxi(t) = 0$).\\

 \item
 If
 \begin{equation}
 \hp(t) = \frac{h_0(t)\hX(t)}{\lambda_0(t)}, \text{i.e.}\; \; \hX(t) = \frac{\lambda_0(t) \hp(t)}{h_0(t)},
 \end{equation}
 then we harvest immediately from $\hX(t)$ at a rate $d\hxi(t)$ which is exactly enough to prevent $\hX(t)$ from going above $\frac{\lambda_0(t) \hp(t)}{h_0(t)}.$\\

 \item
If
 \begin{equation}
 \hp(0^-) < \frac{h_0(0)\hX(0^-)}{\lambda_0(0)},  \text{i.e.}\; \; \hX(0^-) > \frac{\lambda_0(0) \hp(0^-)}{h_0(0)},
 \end{equation}
 then we harvest immediately what is necessary to bring  $\hX(0)$ down to the the level of $\frac{\lambda_0(0)\hp(0}{h_0(0)}.$\\

 \end{itemize}

\end{theorem}

\subsection{Application to model uncertainty singular control}\label{sec6.2}
We  represent {\em model uncertainty}  by a family of probability measures $Q = Q^\theta$  equivalent to $P$, with the Radon-Nikodym derivative on ${\mathcal  F}_t$ given by
\begin{equation}\label{eq6.1}
\frac{d(Q \mid {\mathcal  F}_t)}{d(P \mid {\mathcal  F}_t)} = G^\theta(t)
\end{equation}
where, for $0 \leq t \leq T$, $G^\theta(t)$ is an exponential martingale of the form
\begin{equation}\label{eq6.7}
dG^\theta(t) = G^\theta(t^-) \theta(t)  dB(t);
\;\;
G^\theta(0) = 1.
\end{equation}
Here $\theta$ may be regarded as a {\em scenario control}.
Let $\AC_1:={\AC_\GB}$ denote a given family of admissible singular controls $\xi$ and let $\AC_2:=\Theta$ denote a given set of admissible scenario controls $\theta$
 such that
 \begin{equation}\label{eq6.8}
E [ \int_0^T  \theta^2(t) dt] < \infty.
 \end{equation}


Now assume that $X_1(t) = X^{\xi}(t)$ is a singularly controlled mean-field It\^o process of the form
\begin{align}\label{eq6.10}
dX_1(t) & = b_1(t,X_1(t),Y_1(t),\omega)dt + \sigma_1(t,X_1(t),Y_1(t),\omega)dB(t) \nonumber \\
&+ \lambda_1(t,X_1(t),\omega)d\xi(t),
\end{align}
where
\begin{equation}\label{eq6.11}
Y_1(t) = F(X_1(t, \cdot)),
\end{equation}
and $F$ is a Fr\'echet differentiable operator on $L^2(P)$.

As before let $\GB^{(1)}= \{ \GC_t^{(1)} \}_{0 \leq t \leq T}$ and $\GB^{(2)}=  \{\GC_t^{(2)} \}_{0 \leq t \leq T}$  be  given subfiltrations of $\FB= \{ \FC_t \}_{0 \leq t \leq T}$, representing the information available to the controllers at time $t$. It is required that $\xi \in \AC_1$ be $\GB^{(1)}$-predictable, and  $\theta \in \AC_2$ be $\GB^{(2)}$-predictable.  We set
$w = (\xi,\theta)$ and consider the stochastic differential  game to find $(\hat{\xi}, \hat{\theta}) \in {\mathcal  A}_1 \times {\mathcal  A}_2$ such that
\begin{equation}\label{eq6.12}
\sup_{\xi \in {\mathcal  A}_1}\inf_{\theta \in {\mathcal  A}_2}
E_{Q^\theta}[j(\xi,\theta)] =
E_{Q^{\hat{\theta}}}[j(\hat{\xi},\hat{\theta})] = \inf_{\theta \in {\mathcal  A}_2} \sup_{\xi \in {\mathcal  A}_1}
E_{Q^\theta}[j(\xi,\theta)],
  \end{equation}
where
\begin{align}\label{eq6.13}
j(\xi,\theta)
&=  \int_0^T \{f_1(t,X(t),Y(t),\xi(t),\omega)+\rho(\theta(t)) \}dt \nonumber\\
& \qquad
+ g_1(X(T),Y(T),\omega) + \int_0^T h_1(t,X(t),\omega) \xi(dt).
 \end{align}

The term $E_{Q^\theta} [\int_0^T \rho(\theta(t)) dt]$
can be seen as  a penalty term, penalizing the difference between $Q^\theta$ and the original probability measure $P$.

Note that since $G^\theta(t)$ is a martingale we have
\begin{align}\label{eq6.14}
E_{Q^\theta}[j(\xi,\theta)] & = E \left[ G^\theta(T) g_1(X(T),Y(T)) + \int_0^T G^\theta(t) \{ f_1(t,
X(t), Y(t),\xi(t))) + \rho(\theta(t)) \} dt \right. \nonumber\\
& + \left. \int_0^T G^\theta (t) h_1(t,X(t),Y(t)) \xi(dt)\right] =: J(\xi,\theta).
\end{align}

We see that this is a mean-field singular control stochastic differential game of the type discussed in Section 4, with a two-dimensional state space
\begin{equation}\label{eq6.15}
X(t) := (X_1(t), X_2(t)):= (X^\xi(t),G^\theta(t))
\end{equation}
and with
\begin{align}\label{eq6.16}
f(t,X(t),Y(t),\xi,\theta)) &:= G^\theta(t) \{f_1(t,X_1(t),Y_1(t),\xi(t))+\rho(\theta(t))\} \nonumber\\
& =X_2(t) \{f_1(t,X_1(t),Y_1(t),\xi(t))+\rho(\theta(t))\},
\end{align}

\begin{equation}\label{eq6.17}
 g(X(T),Y(T)):=G^\theta(T) g_1(X_1(T),Y_1(T))=X_2(T) g_1(X_1(T),Y_1(T)),
 \end{equation}
  and
 \begin{equation}\label{6.18}
 h(t,X(t),Y(t)):= G^\theta(t) h_1(t,X_1(t),Y_1(t))=X_2(t) h_1(t,X_1(t),Y_1(t)).
\end{equation}
Using the result from Section 4, we get the following Hamiltonian for the game \eqref{eq6.12}:
\begin{align}\label{eq6.19}
 H& (t,x_1,x_2,y_1,\xi,\theta,p,q,r)(dt, \xi(dt)) \nonumber \\
 & = H_{0} (t,x_1,x_2,y_1,\xi,\theta,p,q) dt + \{ \lambda_1(t,x) p_1 + x_2h_1(t,x)\} \xi(dt) 
 \end{align}
 where
 \begin{align}\label{eq6.20}
H_{0}(t,x_1,x_2,y_1,\xi,\theta,p,q)& =x_2 \{f_1(t,x_1,y_1,\xi) + \rho(\theta) \} \nonumber\\
& + b_1(t,x_1,y_1) p_1 + \sigma(t,x_1,y_1)q_1+ x_2 \theta q_2. 
 \end{align}
The corresponding mean-field BSDEs for the adjoint processes become
\begin{equation}\label{eq6.21}
 \begin{cases}
 dp_1(t) & = \displaystyle - \frac{\partial H_{0}}{\partial x_1} (t,X(t),Y(t),w(t),p(t),q(t))dt\\
 &\displaystyle- \frac{\partial H_{0}}{\partial y_1} (t,X(t), Y(t), w(t),p(t),q(t))\nabla_{X_1(t)} F)dt
  + q_1(t) dB(t) \\
  p_1(T) & = \displaystyle X_2(T) \frac{\partial g_1}{\partial x_1}(X_1(T),Y_1(T))+E[\frac{\partial g_1}{\partial y_1}(X_1(T),Y_1(T))]\nabla_{X_1(T)} F\\
\end{cases}
\end{equation}
and
\begin{equation}\label{eq6.22}
 \begin{cases}
 dp_2(t) & = \displaystyle - \{f_1(t,X_1(t),Y_1(t),\xi) + \rho(\theta(t)) + \theta(t) q_2(t) \}dt
 + q_2(t) dB(t) \\
  p_2(T) & = \displaystyle  g_1(X_1(T),Y_1(T)).\\
\end{cases}
\end{equation}
Minimizing the Hamiltonian with respect to $\theta$ gives the following first order condition:
\begin{equation}\label{eq6.23}
\frac{\partial \rho}{\partial \theta_0}(t) = -E[q_2(t) \mid \GC^{(2)}_t].
\end{equation}
The variational inequalities \eqref{eq3.13a} - \eqref{eq3.13b} reduce to
\begin{align}\label{eq6.25}
\begin{cases}
& \left. E[\nabla_{\xi}f_1(t,\hX(t),\hY(t),\hw(t))+\lambda_1(t,\hX_1(t))\hp(t) + h_1(t,\hX_1(t)) \right.  \left.\mid \EC_t^{(1)} ] \leq 0;  \right.\\
 & \left. E[\nabla_{\xi}f_1(t,\hX(t),\hY(t),\hw(t))+\lambda_1(t,\hX(t))\hp(t) + h_1(t,\hX(t)) \right.  \left. \mid \EC_t^{(1)} ] d\hxi(t) = 0;   \right.
 \end{cases}
 \end{align}

\subsection{A special case}\label{se6.3}

For simplicity, consider the special case with
\begin{equation} \label{eq6.26}
\GC^{(i)}_t = \FC_t, \;\;
\end{equation}
and
\begin{equation} \label{eq6.27}
 \lambda_1(t,x) = \lambda_1(t), \; \; h_1(t,x) = h_1(t)
\end{equation}
i.e., $ \lambda_1$ and $h_1$ do not depend on $x$.

Then, writing $X_1(t) = X^{\xi}(t)$, $Y_1(t) = Y^{\xi}(t)$ and $X_2(t) = G^{\theta}(t)$ and $f_1 = f, g_1=g, b_1=b,\sigma_1=\sigma, \lambda_1=\lambda$,
 the controlled system gets the form $(X^{\xi},G^{\theta})$, where $G^{\theta}$ is given by \eqref{eq6.7} and 
\begin{equation}\label{eq6.28}
dX^{\xi}(t)  = b(t,X^{\xi}(t),Y^{\xi}(t))dt + \sigma(t,X^{\xi}(t),Y^{\xi}(t))dB(t) + \lambda(t)d\xi(t); \;\;X^{\xi}(0)=x  \\
\end{equation}
The performance functional becomes
\begin{align} \label{eq6.30}
E_{Q^\theta}[j(\xi,\theta)] & = E \left[ G^\theta(T) g(X^\xi(T),Y\xi(T))\right. \nonumber\\
&\left. \qquad + \int_0^T G^\theta(t) \{ f(t,
X^\xi(t), Y^\xi(t),\xi(t)) + \rho(\theta(t)) \} dt] =: J(\xi,\theta).\right.
\end{align}
and the Hamiltonian becomes
\begin{align}
 H& (t,x,g,y,\xi,\theta,p,q)(dt, \xi(dt)) \nonumber \\
 & = H_{0} (t,x,g,y,\xi,\theta,p,q) dt + \{ \lambda(t) p_1 + gh_1(t)\} \xi(dt), \label{eq6.31}
  \end{align}
 where
 \begin{equation}\label{eq6.32}
H_{0}(t,x,g,y,\xi,\theta,p,q) =g \{f_1(t,x,y,\xi) + \rho(\theta) \} + b_1(t,x,y) p_1 + \sigma(t,x,y)q_1+ g \theta q_2. \nonumber\\
 \end{equation}
The corresponding mean-field BSDEs for the adjoint processes become
\begin{equation}\label{eq6.33}
 \begin{cases}
 dp_1(t) & = \displaystyle - \frac{\partial H_{0}}{\partial x} (t,X^\xi(t),Y^\xi(t),\xi(t),p(t),q(t))dt\\
 &\displaystyle- \frac{\partial H_{0}}{\partial y} (t,X^\xi(t), Y^\xi(t), \xi(t),p(t),q(t))\nabla_{X^\xi(t)} F)dt+ q_1(t) dB(t) \\
 p_1(T) & = \displaystyle G^\theta(T) \frac{\partial g}{\partial x}(X^\xi(T),Y^\xi(T))+E[\frac{\partial g}{\partial y}(X^\xi(T),Y^\xi(T))]\nabla_{X^\xi(T)} F\\
\end{cases}
\end{equation}
and
\begin{equation}\label{eq6.34}
 \begin{cases}
 dp_2(t) & = \displaystyle - \{f(t,X^{\xi}(t),Y^{\xi}(t),\xi) + \rho(\theta(t)) + \theta(t) q_2(t) \}dt+ q_2(t) dB(t)\\

  p_2(T) & = \displaystyle  g_1(X^{\xi}(T),Y^{\xi}(T)).\\
\end{cases}
\end{equation}

Then the first order condition for a minimum of the Hamiltonian with respect to $\theta_0$ reduces to
\begin{equation}\label{eq6.35}
\rho'(\theta)(t) = -q_2(t).
\end{equation}

The variational inequalities \eqref{eq6.19} become
\begin{align}\label{eq6.36}
\begin{cases}
& \left. \nabla_{\xi}f_1(t,X(t),Y(t),\xi(t))+\lambda_1(t)p_1(t) + h_1(t) \leq 0;  \right.\\
 & \left. [\nabla_{\xi}f_1(t,X_1(t),Y_1(t),\xi(t))+\lambda_1(t)p_1(t) + h_1(t)] d\xi(t) = 0;   \right.
 \end{cases}
 \end{align}

In general it seems to be a formidable mathematical challenge to solve such a coupled system of forward-backward singular SDEs. However, in some cases a possible solution procedure could be described, as in the next example:

\subsection{Optimal harvesting under uncertainty}\label{sec6.4}

Now we consider a model uncertainty version of the optimal harvesting problem in Section 6.1. For simplicity we put $K=1$. Thus we have the following mean-field forward system $(X,^{\xi},G^{\theta})$, where $G^{\theta}$ is given by \eqref{eq6.7} and 
\begin{equation}\label{eq6.37}
\begin{cases}
dX(t) &= dX^\xi(t)= E[X^\xi(t)]b(t)dt + X^\xi(t) \sigma(t)dB(t) - \lambda_0(t) d\xi(t)\\
X^\xi(0^-)&=x > 0, \\
\end{cases}
\end{equation}

with performance functional
\begin{equation} \label{eq6.39}
J(\xi,\theta) = E[\int_0^T G^{\theta}(t) \rho(\theta(t))dt + \int_0^T G^{\theta}(t) h_0(t)X^\xi(t) d\xi(t) + G^{\theta}(T)X^\xi(T)].
\end{equation}

The model uncertainty harvesting problem is to find $(\xi^*, \theta^*) \in \AC_\FC \times \Theta$ such that
\begin{equation}\label{eq6.39a}
\sup_{\xi \in \AC_\FC} \left\{ \inf_{\theta \in \Theta} J(\xi, \theta)\right\} = \inf_{\theta \in \Theta} \left\{ \sup_{\xi \in \AC_\FC} J(\xi, \theta)\right\} = J (\xi^*, \theta^*),
\end{equation}

Here the Hamiltonian is
\begin{align}\label{eq6.40}
&H(t,x,g,y,\xi,\theta,p,q)(dt, \xi(dt)) \nonumber\\
&= \{g \rho(\theta) + yb(t) p_1 + x\sigma(t)q_1+ g \theta q_2\} dt + \{ -\lambda_0(t) p_1 + xgh_0(t)\} \xi(dt).
 \end{align}

Minimizing the Hamiltonian with respect to $\theta$ gives the first order equation
\begin{equation} \label{eq6.41}
\rho'(\theta)(t) = -q_2(t).
\end{equation}

The corresponding reflected backward system is
\begin{equation}\label{eq6.42}
 \begin{cases}
 dp_1(t) & = \displaystyle - [b(t)p_1(t) + \sigma(t)q_1(t)]dt -h_0(t)G^\theta(t) d\xi(t) + q_1(t) dB(t)\\
 p_1(T) & = \displaystyle  G^\theta(T) \\
\end{cases}
\end{equation}
and
\begin{equation}\label{eq6.43}
 \begin{cases}
 dp_2(t) & = \displaystyle - [ \rho(\theta(t)) + \theta(t) q_2(t)] dt-h_0(t)X(t) d\xi(t)+ q_2(t) dB(t)\\
  p_2(T) & = \displaystyle  X(T),\\
\end{cases}
\end{equation}
with variational inequalities
\begin{equation}\label{eq6.44}
\begin{cases}
&  -\lambda_0(t)p_1(t) + h_1(t) \leq 0;  \\
 & [ -\lambda_0(t)p_1(t) + h_1(t)] d\xi(t) = 0
 \end{cases}
 \end{equation}
Then we get the following result:

\begin{theorem}
Suppose there exists a solution $\hX(t):=X^{\hxi}(t)$, $\hG(t):= G^{\hth}(t)$, $\hp_1(t)$,$\hq_1(t)$,$\hp_2(t)$,$\hq_2(t)$,$\hxi(t)$,$\hth(t)$ of the coupled system of mean-field forward-backward singular stochastic differential equations consisting of the forward equations \eqref{eq6.37}
and the reflected backward equations \eqref{eq6.42},\eqref{eq6.43}, and satisfying the constraint \eqref{eq6.44}. 
Then $\hxi(t)$ is the optimal harvesting strategy and $\hth(t)$ is the optimal scenario parameter for the model uncertainty harvesting problem \eqref{eq6.39a}.
\end{theorem}

\subsection{A mean field singular game}

We now return to the mean field singular game described in Section 2.3.
\vskip 0.3cm

In this case, we get from \eref{eq3.6}
\begin{eqnarray*}
&&H_i(t,x,y,\xi_1, \xi_2, p_i, q_i)(dt, \xi_1(dy), \xi_2(dt))\\
&&\qquad =\pi \min(x, \xi_1+\xi_2)+yb(t)p_i+x\si(t) q_i+h_1\xi_1(dt)+h_2\xi_2(dt)
\end{eqnarray*}
and the adjoint equations \eref{eq3.8} becomes
\[
\left\{ \begin{array}{lll}
dp_i(t)&=&-\left[ \chi_{[0, \xi_1+\xi_2)}(x) \pi(t) +\si(t)q_i(t)+b(t)p_i(t)\right] dt+q_i(t)dB(t)\,;\\
p_i(T)&=&0\,.
\end{array}\right.
\]
The variational inequalities  \eref{eq5.9} get the form
\[
\left\{ \begin{array}{ll }
 \pi(t)\chi_{[0, X(t))}(\xi_1(t)+\xi_2(t))+h_i(t) & \le 0\\
\hbox{and}& \\
 \left[\pi(t)\chi_{[0, X(t))}(\xi_1(t)+\xi_2(t))+h_i(t)\right]\xi_i(dt) &= 0\,;\qquad i=1,2\,.\\
\end{array}\right.
\]
Optimal strategy for factory 1:
\begin{enumerate}
\item[i)]\ If $\pi(t)+h_1(t)<0$, do nothing.
\item[ii)]\ If $\pi(t)+h_1(t)\ge 0$ increase $\xi_1(t)$ to $X(t)-\xi_2^*(t)$.
\end{enumerate}
In other words, $(\xi_1, X)$ solves  the reflected Skorohod problem
\begin{equation}
\left\{\begin{array}{ll}
\xi_1(t)\ge(X(t)-\xi_2^*(t))\chi_{[0, \infty)}(\pi(t)+h_1(t))\\
\left[\xi_1(t)-(X(t)-\xi_2^*(t))\chi_{[0, \infty)}(\pi(t)+h_1(t))\right]\xi_1(dt)=0\,.
\end{array}\right.\label{e.6.45}
\end{equation}
So for given $\xi_2^*$ we choose $\xi_1:=R_1(\xi_2^*)$ solution of the
reflected Skorohod problem \eref{e.6.45}.

Similarly, given $\xi_1^*$ we choose  $\xi_2:=R_2(\xi_1^*)$ as the solution of the
reflected Skorohod problem
\begin{equation}
\left\{\begin{array}{ll}
\xi_2(t)\ge(X(t)-\xi_1^*(t))\chi_{[0, \infty)}(\pi(t)+h_2(t))\\
\left[\xi_2(t)-(X(t)-\xi_1^*(t))\chi_{[0, \infty)}(\pi(t)+h_2(t))\right]\xi_2(dt)=0\,.
\end{array}\right.\label{e.6.46}
\end{equation}
Thus, to find the {\it Nash equilibrium} we need to solve the following coupled reflected
Skorohod problem:
\begin{eqnarray}
\left\{\begin{array}{ll}
\xi_1(t)\ge(X(t)-\xi_2 (t))\chi_{[0, \infty)}(\pi(t)+h_1(t))\\
\left[\xi_1(t)-(X(t)-\xi_2 (t))\chi_{[0, \infty)}(\pi(t)+h_1(t))\right]\xi_1(dt)=0\,; \\
\\
\xi_2(t)\ge(X(t)-\xi_1 (t))\chi_{[0, \infty)}(\pi(t)+h_2(t))\\
\left[\xi_2(t)-(X(t)-\xi_1 (t))\chi_{[0, \infty)}(\pi(t)+h_2(t))\right]\xi_2(dt)=0\,.
\end{array}\right. \label{e.6.47}
\end{eqnarray}

The above system of reflected Skorohod problem can be solved in the following way:

We divide the interval$[0, T]$ into $0=t_0<t_1<\cdots<t_n=T $ such that
on each interval $[t_k, t_{k+1})$ the signs of $\pi(t)+h_1(t)$ and $\pi(t)+h_2(t)$ remains
unchanged.

On each interval $[t_k, t_{k+1}]$, we use the following control principles.
If both of the inequalities $\pi(t)+h_1(t)<0$ and  $\pi(t)+h_2(t)<0$ hold,
then do nothing.  If $\pi(t)+h_1(t)<0$ and  $\pi(t)+h_2(t)\ge 0$, then
the first factory does not do anything. The second condition in
 \eref{e.6.47}
becomes
\begin{eqnarray*}
\left\{\begin{array}{ll}
\xi_2(t)\ge (X(t)-\xi_1 (t))  \\
\left[\xi_2(t)-(X(t)-\xi_1 (t))  )\right]\xi_2(dt)=0\,.
\end{array}\right.
\end{eqnarray*}
By Remark 2.7 (namely, Equation (2.8)) of \cite{bkr},
\[
\xi_2(t)=  \sup_{t_k\le s\le t}(X(s)-\xi_1 (t_k))^+ \,,\quad t_k\le t\le t_{k+1}\,.
\]
Thus, we  keep $\xi_1(t)=\xi_1(t_k)$ unchanged and
in the same time   increase $\xi_2(t)$ to $X(t)-\xi_1 (t)$.
Similar result holds if $\pi(t)+h_1(t)\ge 0$ and  $\pi(t)+h_2(t)< 0$.

If both $\pi(t)+h_1(t)\ge 0$ and $\pi(t)+h_2(t)\ge 0$, then \eref{e.6.47}
becomes
\begin{eqnarray}
\left\{\begin{array}{ll}
\xi_1(t)\ge(X(t)-\xi_2 (t)) \\
\left[\xi_1(t)-(X(t)-\xi_2 (t)) \right]\xi_1(dt)=0\,; \\
\\
\xi_2(t)\ge(X(t)-\xi_1 (t)) \\
\left[\xi_2(t)-(X(t)-\xi_1 (t)) )\right]\xi_2(dt)=0\,.
\end{array}\right. \label{e.6.49}
\end{eqnarray}
Let $\xi(t)=\xi_1(t)+\xi_2(t)$ and  then \eref{e.6.49} is equivalent to
\[
\left\{\begin{array}{ll}
\xi(t)\ge X(t)\\
\left[\xi(t)-  X(t)   \right]\xi (dt) =0\,.
\end{array}\right.
\]
Again by Remark 2.7 (namely, Equation (2.8)) of \cite{bkr},  we have
\[
\xi(t)=   \sup_{t_k\le s\le t} X(s) ^+ \,,\quad t_k\le t\le t_{k+1}\,.
\]
Now we show that any decomposition of $\xi(t)$ into the sum of two nondecreasing processes
$\xi_1(t)$ and $\xi_2(t)$ will solve \eref{e.6.49}.  In fact, assume  $\xi(t)=\xi_1(t)+\xi_2(t)$,
where $\xi_1$ and $\xi_2$ are two nondecreasing processes.  Since $\xi(t)\ge X(t)$
and $\xi_1$ and $\xi_2$ are nondecreasing, we have
\[
\left\{\begin{array}{ll}
\left[\xi(t)-  X(t)   \right]\xi_1 (dt)\ge &0\\
\left[\xi(t)-  X(t)   \right]\xi_2 (dt)\ge &0\,.
\end{array}\right.
\]
Add them we have
\[
\left[\xi(t)-  X(t)   \right]\xi_1 (dt)+
\left[\xi(t)-  X(t)   \right]\xi_2 (dt)=\left[\xi(t)-  X(t)   \right]\xi (dt)=0\,.
\]
This implies
\[
\left\{\begin{array}{ll}
\left[\xi(t)-  X(t)   \right]\xi_1 (dt)= &0\\
\left[\xi(t)-  X(t)   \right]\xi_2 (dt)= &0\,.
\end{array}\right.
\]
Thus, $\xi_1$ and $\xi_2$ satisfies \eref{e.6.49}.
Summarizing we have
\begin{theorem} Assume that we can divide the interval $[0, T]$ into $0=t_0<t_1<\cdots<t_n=T $ such that
on each interval $[t_k, t_{k+1})$ the signs of $\pi(t)+h_1(t)$ and $\pi(t)+h_2(t)$ remain
unchanged.  Then we can recursively find the solution $\xi_1$ and $\xi_2$ on each interval $[t_k, t_{k+1}]$
for $k=0, 1, \cdots, n-1$. On the interval $[t_k,  t_{k+1}]$, we have
\begin{enumerate}
\item[(i)] If both of the inequalities $\pi(t)+h_1(t)<0$ and  $\pi(t)+h_2(t)<0$ hold,
then do nothing.
\item[(ii)] If  $\pi(t)+h_1(t)<0$ but $\pi(t)+h_2(t)\ge 0$,  then
\[
\xi_1(t)=\xi_1(t_k)\,, \quad t_k\le t\le t_{k+1}
\]
and
\[
\xi_2(t)=  \sup_{t_k\le s\le t}(X(s)-\xi_1 (t_k))^+ \,,\quad t_k\le t\le t_{k+1}\,.
\]

If  $\pi(t)+h_1(t)\ge 0$ but $\pi(t)+h_2(t)< 0$,  then
\[
\xi_2(t)=\xi_2(t_k)\,, \quad t_k\le  t \le t_{k+1}
\]
and
\[
\xi_1(t)=  \sup_{t_k\le s\le t}(X(s)-\xi_2 (t_k))^+ \,,\quad t_k\le t\le t_{k+1}\,.
\]
\item[(iii)] If  both of the inequalities $\pi(t)+h_1(t)\ge 0$ and  $\pi(t)+h_2(t)\ge 0$ hold,
then   $\xi_1$ and $\xi_2$ can be any nondecreasing processes such that
\[
\xi_1(t)+\xi_2(t)=   \sup_{t_k\le s\le t} X(s) ^+ \,,\quad t_k\le t\le t_{k+1}\,.
\]
\end{enumerate}
Note in particular that in this case the Nash equilibrium is not unique.
\end{theorem}

\end{document}